\documentclass[11pt]{article}
\newcommand{\comment}[1]		{}
\usepackage{amsmath,amsthm}
\usepackage{rotating,pdflscape,afterpage}
\usepackage[sc,osf]{mathpazo}
	\AtBeginDocument{
		\DeclareSymbolFont{AMSb}{U}{msb}{m}{n}
		\DeclareSymbolFontAlphabet{\mathbb}{AMSb}
	}
\usepackage{graphicx}
\usepackage{subcaption}
\usepackage{extpfeil}
\usepackage{mathabx} 
\usepackage{fancyvrb}
\usepackage{float}
\usepackage{caption}
\newcommand{\mockalph}[1]{\!}
\usepackage{enumitem}
\usepackage[nouppercase]{scrlayer-scrpage}
\usepackage{avant} 
\usepackage{mathrsfs}
\usepackage{mathtools}
\usepackage{faktor,xfrac}
\usepackage[all]{xypic}\xyoption{rotate}
	\newdir{ >}{{}*!/-7pt/\dir{>}}
	\newdir{i}{{}*!/-7pt/\dir{(}	}

\usepackage{thmtools}

\usepackage[utf8]{inputenc}
\usepackage[T1]{fontenc}
\usepackage{tocloft}
\makeatletter
\renewcommand{\l@figure}{\@dottedtocline{1}{1em}{3.5em}}
\renewcommand{\l@table}{\@dottedtocline{2}{1em}{3.5em}}
\makeatother
\newcommand*{\noaddvspace}{\renewcommand*{\addvspace}[1]{}}
\addtocontents{lof}{\protect\noaddvspace}

\usepackage[hyphens]{url}		%
\usepackage{hyperref}	
\usepackage{cleveref} 	

\makeatletter
\let\c@figure\c@table
\let\c@equation\c@table
\makeatother

\numberwithin{table}{section}
\numberwithin{figure}{section}
\newtheorem{abcthm}[table]{Theorem}

\newtheorem{theorem}[table]{Theorem}

\newtheorem{proposition}[table]{Proposition}

\newtheorem{corollary}[table]{Corollary}

\newtheorem{lemma}[table]{Lemma}

\newtheorem{claim}[table]{Claim}

\theoremstyle{definition}

\newtheorem{definition}[table]{Definition}

\newtheorem{notation}[table]{Notation}

\newtheorem{observation}[table]{Observation}

\newtheorem{conjecture}[table]{Conjecture}

\newtheorem{discussionsilent}[table]{}

\theoremstyle{remark}
\newtheorem{fact}[table]{Fact}

\newtheorem{example}[table]{Example}
\newtheorem{exercise}[table]{Exercise}

\newtheorem{problem}[table]{Problem}
\newtheorem{histrmks}[table]{Historical remarks}
\newtheorem{remark}[table]{Remark}
\newtheorem{remarks}[table]{Remarks}

\theoremstyle{plain}
\newtheorem*{thm*}{Theorem}
\newtheorem*{theorem*}{Theorem}
\newtheorem*{prop*}{Proposition}
\newtheorem*{proposition*}{Proposition}
\newtheorem*{lemma*}{Lemma}
\newtheorem*{corollary*}{Corollary}
\newtheorem*{cor*}{Corollary}

\theoremstyle{definition}
\newtheorem*{definition*}{Definition}
\newtheorem*{defn*}{Definition}
\newtheorem*{QQ*}{Question}
\newtheorem*{obs*}{Observation}
\newtheorem*{notation*}{Notation}
\newtheorem*{discussion*}{Discussion}

\theoremstyle{remark}
\newtheorem*{rmk*}{Remark}
\newtheorem*{remark*}{Remark}

\newtheorem*{examples*}{Examples}
\newtheorem*{example*}{Example}
\newtheorem*{EG*}{Example}
\newtheorem*{EGs*}{Examples}
\newtheorem*{fact*}{Fact}
\newtheorem*{prob*}{Problem}

\usepackage{etoolbox}
\usepackage[dvipsnames,svgnames,x11names,table]{xcolor}
\newcommand		{\defd}[1]	{\textcolor{RoyalBlue}{\textbf{\textit{#1}}}}
\newcommand		{\defm}[1]	{\textcolor{RoyalBlue}{#1}}

\tracingpatches
\makeatletter
\patchcmd{\@setref}{\bfseries ??}{\bfseries\color{red} FIX ME!}{}{}
\patchcmd{\@setcite}{\bfseries ?}{\bfseries\color{red} FIX ME!}{}{}
\patchcmd{\@setcref}         {??}{\color{red} FIX ME!}{}{}
\patchcmd{\@setcref}         {??}{\color{red} FIX ME!}{}{}
\patchcmd{\@setcrefrange}    {??}{\color{red} FIX ME!}{}{}
\patchcmd{\@setcrefrange}    {??}{\color{red} FIX ME!}{}{}
\patchcmd{\@setcrefrange}    {??}{\color{red} FIX ME!}{}{}
\patchcmd{\@setcrefrange}    {??}{\color{red} FIX ME!}{}{}
\patchcmd{\@setcrefrange}    {??}{\color{red} FIX ME!}{}{}
\patchcmd{\@setcrefrange}    {??}{\color{red} FIX ME!}{}{}
\patchcmd{\@setnamecref}     {??}{\color{red} FIX ME!}{}{}
\patchcmd{\@setnamecref}     {??}{\color{red} FIX ME!}{}{}
\patchcmd{\@setcpageref}     {??}{\color{red} FIX ME!}{}{}
\patchcmd{\@setcpageref}     {??}{\color{red} FIX ME!}{}{}
\patchcmd{\@setcpagerefrange}{??}{\color{red} FIX ME!}{}{}
\patchcmd{\@setcpagerefrange}{??}{\color{red} FIX ME!}{}{}
\patchcmd{\@setcpagerefrange}{??}{\color{red} FIX ME!}{}{}
\patchcmd{\@setcpagerefrange}{??}{\color{red} FIX ME!}{}{}
\patchcmd{\@setcpagerefrange}{??}{\color{red} FIX ME!}{}{}
\patchcmd{\@cref}            {??}{\color{red} FIX ME!}{}{}
\def\blx@citation@entry#1#2{%
  \blx@bibreq{#1}%
  \ifinlist{#1}{\blx@cites}
    {}
    {\listgadd{\blx@cites}{#1}%
     \blx@auxwrite\@mainaux{}{\string\abx@aux@cite{#1}}}%
  \ifinlistcs{#1}{blx@segm@\the\c@refsection @\the\c@refsegment}
    {}
    {\listcsgadd{blx@segm@\the\c@refsection @\the\c@refsegment}{#1}}%
  \blx@ifdata{#1}%
    {}%
    {\ifcsdef{blx@miss@\the\c@refsection}%
       {\ifinlistcs{#1}{blx@miss@\the\c@refsection}%
          {{\bfseries\color{red} cite:} }%
          {\blx@logreq@active{#2{#1}}}}%
       {\blx@logreq@active{#2{#1}}}}}
\def\blx@citeadd#1{%
  \ifcsdef{blx@keyalias@\the\c@refsection @#1}
    {\edef\blx@realkey{\csuse{blx@keyalias@\the\c@refsection @#1}}}
    {\def\blx@realkey{#1}}%
  \expandafter\blx@citation\expandafter{\blx@realkey}\blx@msg@cundefon
  \expandafter\blx@ifdata\expandafter{\blx@realkey}
    {\advance\blx@tempcnta\@ne
     \listeadd\blx@tempa{\blx@realkey}}
    {\ifnum\blx@tempcntb>\z@\multicitedelim\fi
     \expandafter\abx@missing\expandafter{\blx@realkey}%
     \advance\blx@tempcntb\@ne}}

\DeclarePairedDelimiterX{\pmodx}[1]{(}{)}{{\operator@font mod}\mkern6mu#1}
\renewcommand{\pmod}{%
  \allowbreak
  \if@display\mkern18mu\else\mkern8mu\fi
  \pmodx
}
\makeatother

\makeatletter
\newcommand{\oset}[3][0ex]{%
\raisebox{.175ex}{$%
  \mathrel{\mathop{#3}\limits^{
    \vbox to#1{\kern-2\ex@
    \hbox{$\scriptstyle#2$}\vss}}}
    $}%
    }
\makeatother

\newcommand{\myred}{BrickRed}
\hypersetup{
		colorlinks			= true, 	
		urlcolor			= \myred, 
		linkcolor			= \myred, 	
		citecolor			= PineGreen 	
}

\usepackage{xspace}

\usepackage[left=2.54cm,right=2.54cm,top=2.54cm,bottom=2.54cm]{geometry}
\usepackage{tikz}
\usetikzlibrary{matrix,fit,backgrounds,calc,arrows} 

\tikzstyle{image}=[rectangle,fill=Red!20,inner sep=-2pt]
\tikzstyle{nonzero}=[rectangle,fill=Navy!20,inner sep=0pt]
\tikzstyle{nonzerosm}=[rectangle,fill=Navy!20,inner sep=-2pt]

\makeatletter
\newbox\xrat@below
\newbox\xrat@above
\newcommand{\xrightarrowtail}[2][]{%
  \setbox\xrat@below=\hbox{\ensuremath{\scriptstyle #1}}%
  \setbox\xrat@above=\hbox{\ensuremath{\scriptstyle #2}}%
  \pgfmathsetlengthmacro{\xrat@len}{max(\wd\xrat@below,\wd\xrat@above)+.6em}%
  \mathrel{\tikz [>->,baseline=-.55ex]
                 \draw (0,0) -- node[below=-2pt] {\box\xrat@below}
                                node[above=-2pt] {\box\xrat@above}
                       (\xrat@len,0) ;}}
\newbox\xrat@below
\newbox\xrat@above
\renewcommand{\xtwoheadrightarrow}[2][]{%
  \setbox\xrat@below=\hbox{\ensuremath{\scriptstyle #1}}%
  \setbox\xrat@above=\hbox{\ensuremath{\scriptstyle #2}}%
  \pgfmathsetlengthmacro{\xrat@len}{max(\wd\xrat@below,\wd\xrat@above)+.6em}%
  \mathrel{\tikz [->>,baseline=-.55ex]
                 \draw (0,0) -- node[below=-2pt] {\box\xrat@below}
                                node[above=-2pt] {\box\xrat@above}
                       (\xrat@len,0) ;}}
\makeatother
\newcommand{\xmono}{\xrightarrowtail}
\newcommand{\mono}{\xmono{\phantom{\ \, }}}
\newcommand{\xepi}{\xtwoheadrightarrow}
\newcommand{\epi}{\xepi{\phantom{\ \, }}}

\makeatletter
\newcommand{\presectionskip}{-1.5\baselineskip}
\newcommand{\postsectionskip}{0.3\baselineskip}
\usepackage{titlesec}
\renewcommand{\section}{\@startsection
  {chapter}{0}{0mm}
  {\presectionskip}
  {\postsectionskip}
  {\sffamily\huge}}
\renewcommand{\section}{\@startsection
  {section}{1}{0mm}
  {\presectionskip}
  {\postsectionskip}
  {\sffamily\LARGE}}
\renewcommand{\subsection}{\@startsection
  {subsection}{2}{0mm}
  {\presectionskip}
  {\postsectionskip}
  {\sffamily\Large}}
\renewcommand{\subsubsection}{\@startsection
  {subsubsection}{3}{0mm}
  {\presectionskip}
  {\postsectionskip}
  {\sffamily\normalsize}}
\renewcommand{\@seccntformat}[1]{\csname the#1\endcsname.\quad}
\newcommand\HUGE{\@setfontsize\Huge{30}{47}} 
\titleformat{\chapter}[display]
{\sffamily\Large}
{Chapter {\HUGE\normalfont\thechapter}}    
{1em}
{\huge}
\titleformat{\section}
{\sffamily\Large}
{\normalfont{\@setfontsize\Large{18}{47}\thesection.}}    
{1em}
{}
\titleformat{\subsection}
{\sffamily\large}
{\normalfont{\Large\thesubsection}.}    
{1em}
{}
\titleformat{\subsubsection}
{\sffamily\normalsize}
{\normalfont{\large\thesubsubsection}.}    
{1em}
{}
\makeatother

\makeatletter
\def\smallunderbrace#1{\mathop{\vtop{\m@th\ialign{##\crcr
   $\hfil\displaystyle{#1}\hfil$\crcr
   \noalign{\kern3\p@\nointerlineskip}%
   \tiny\upbracefill\crcr\noalign{\kern3\p@}}}}\limits}
\makeatother

\ExplSyntaxOn
\NewDocumentEnvironment{adjunctions}{O{}}
{
	\cs_set_eq:cN {@arraycr} \farin_arraycr:
	\keys_set:nn { farin/adjunction } { #1 }
	\begin{array}
		{
			@{ \hspace { \dim_eval:n { \l_farin_left_shift_dim + \l_farin_padding_dim } } }
			r
			@{ {\farin_strut:} \l_farin_symbol_tl {} }
			l
			@{ \hspace { \dim_eval:n { \l_farin_right_shift_dim + \l_farin_padding_dim } } }
		}
	}
	{
	\end{array}
}
\keys_define:nn { farin/adjunction }
{
	leftshift       .dim_set:N = \l_farin_left_shift_dim,
	leftshift       .initial:n = 0pt,
	rightshift      .dim_set:N = \l_farin_right_shift_dim,
	rightshift      .initial:n = 0pt,
	padding         .dim_set:N = \l_farin_padding_dim,
	padding         .initial:n = 6pt,
	symbol          .tl_set:N  = \l_farin_symbol_tl,
	symbol          .initial:n = \longrightarrow,
	verticalspacing .dim_set:N  = \l_farin_vertspac_dim,
	verticalspacing .initial:n = {3pt},
}
\cs_new_protected:Npn \farin_strut:
{
	\vrule height \dim_eval:n { \ht\strutbox + 1.2\l_farin_vertspac_dim }
	depth  \dim_eval:n { \dp\strutbox + \l_farin_vertspac_dim }
	width 0pt
}
\makeatletter
\exp_args:NNo \cs_new:Npn \farin_arraycr:
{
	\@arraycr\hline
}
\makeatother
\ExplSyntaxOff

\renewcommand{\SS}{\textsection}
\newcommand{\bthm}{\begin{theorem}}
\newcommand{\ethm}{\end{theorem}}
\newcommand{\bprop}{\begin{proposition}}
\newcommand{\eprop}{\end{proposition}}
\newcommand{\bcor}{\begin{corollary}}
\newcommand{\ecor}{\end{corollary}}
\newcommand{\bconj}{\begin{conjecture}}
\newcommand{\econj}{\end{conjecture}}
\newcommand{\blem}{\begin{lemma}}
\newcommand{\elem}{\end{lemma}}
\newcommand{\bclm}{\begin{claim}}
\newcommand{\eclm}{\end{claim}}
\newcommand{\bpf}{\begin{proof}}
\newcommand{\epf}{\end{proof}}
\newcommand{\bdetails}{\begin{details}}
\newcommand{\edetails}{\end{details}}
\newcommand{\bdefi}{\begin{definition}}
\newcommand{\edefi}{\end{definition}}
\newcommand{\bdefn}{\begin{definition}}
\newcommand{\edefn}{\end{definition}}
\newcommand{\bex}{\begin{example}}
\newcommand{\eex}{\end{example}}
\newcommand{\bprob}{\begin{problem}}
\newcommand{\eprob}{\end{problem}}
\newcommand{\bob}{\begin{observation}}
\newcommand{\eob}{\end{observation}}
\newcommand{\bexer}{\begin{exercise}}
\newcommand{\eexer}{\end{exercise}}
\newcommand{\bexers}{\begin{exercises}}
\newcommand{\eexers}{\end{exercises}}
\newcommand{\brmk}{\begin{remark}}
\newcommand{\ermk}{\end{remark}}
\newcommand{\bhist}{\begin{histrmks}}
\newcommand{\ehist}{\end{histrmks}}
\newcommand{\brmks}{\begin{remarks}}
\newcommand{\ermks}{\end{remarks}}
\newcommand{\bverb}{\begin{verbatim}}
\newcommand{\everb}{\end{verbatim}}
\newcommand{\bntn}{\begin{notation}}
\newcommand{\entn}{\end{notation}}
\newcommand{\bfct}{\begin{fact}}
\newcommand{\efct}{\end{fact}}
\newcommand{\bfcts}{\begin{facts}}
\newcommand{\befcts}{\end{facts}}
\newcommand{\benum}{\begin{enumerate}}
\newcommand{\eenum}{\end{enumerate}}
\newcommand{\bitem}{\begin{itemize}}
\newcommand{\eitem}{\end{itemize}}

\renewcommand	{\o}		{\circ}

\renewcommand	{\epsilon}	{\varepsilon}
\renewcommand	{\a}		{\alpha}
\renewcommand	{\b}		{\beta}
\renewcommand	{\d}		{\delta}

\renewcommand	{\l}		{\lambda}

\newcommand		{\vk}		{\varkappa}

\renewcommand	{\:}		{\colon}

\newcommand		{\quotientmedddd}[2]	{{\raisebox{.05em}{$#1$}}\  \!\!\big/\!\!\ 
	{\raisebox{-.05em}{$#2$}}}

\newcommand		{\qquotientmed}[2]	{{\raisebox{.2em}{$#1$}} \ \big/\!\!\!\!\!\big/ \ 
										{\raisebox{-.2em}{$#2$}}}
\newcommand		{\qquotientmedd}[2]	{{\raisebox{.1em}{$#1$}} \ \big/\!\!\!\!\!\big/ \ 
	{\raisebox{-.1em}{$#2$}}}

\newcommand		{\fs}		{{\mathfrak s}}

\newcommand		{\wP}		{\widehat{P}}

\newcommand		{\ewP}	{\ext\wP}

\newcommand		{\SSS}	{Serre spectral sequence\xspace}
\newcommand		{\EMSS}	{Eilenberg--Moore spectral sequence\xspace}

\newcommand		{\eqfity}	{equivariant formality\xspace}
\newcommand		{\isotf}	{isotropy-formal\xspace}
\newcommand		{\isotfity}	{isotropy-formality\xspace}

\newcommand		{\exterior}	{\Lambda}
\newcommand		{\ext}		{\exterior}

\newcommand		{\susp}		{\Sigma}

\newcommand	{\kk}		{k}

\newcommand		{\E}		{E^*}

\newcommand		{\ab}			{^{\mathrm{ab}}}

\newcommand		{\CGA}		{\textsc{cga}\xspace}

\newcommand		{\DGA}		{\textsc{dga}\xspace}

\newcommand		{\CDGA}		{\textsc{cdga}\xspace}

\newcommand	{\fa}{\f a}
\newcommand	{\fb}{\f b}

\newcommand{\Hp}{H^{**}}

\newcommand{\wR}{\wh{R}}

\newcommand{\compl}{\!\wh{\ ^{\phantom{x}}}\mn}
\newcommand{\K}{K^*}
\newcommand{\KG}{\K_G}
\newcommand{\KT}{\K_T}
\newcommand{\KS}{\K_S}

\DeclareMathOperator{\Gr}{Gr}

\newcommand{\WEF}{weakly equivariantly formal\xspace}

\newcommand{\ccpair}{compact, connected pair\xspace}

\newcommand{\GK}{$\smash{(G,K)}$\xspace}
\newcommand{\GS}{$\smash{(G,S)}$\xspace}

\newcommand{\CIR}{complete intersection ring\xspace}
\newcommand{\TFAE}{The following are equivalent\xspace}

\DeclareMathOperator{\gr}		{gr }

\newcommand		{\blbl}		{^{\bl,\bl}}

\newcommand		{\Top}	{\mathsf{Top}}

\makeatletter
\newcommand{\subalign}[1]{%
  \vcenter{%
    \Let@ \restore@math@cr \default@tag
    \baselineskip\fontdimen10 \scriptfont\tw@
    \advance\baselineskip\fontdimen12 \scriptfont\tw@
    \lineskip\thr@@\fontdimen8 \scriptfont\thr@@
    \lineskiplimit\lineskip
    \ialign{\hfil$\m@th\scriptstyle##$&$\m@th\scriptstyle{}##$\crcr
      #1\crcr
    }%
  }
}
\makeatother

\newcommand		{\bl}		{\bullet}

\newcommand		{\Ei}		{E_\infty}
\newcommand		{\col}		{^{0,\bl}}

\newcommand		{\ang}[1]			{\langle #1 \rangle}

\newcommand		{\quation}[1]		{\begin{equation} #1 \end{equation}}

\newcommand		{\eqn}[1]			{\begin{align*} #1 \end{align*}}

\newcommand		{\hyref}[1]			{\hyperref[#1]{\ref*{#1}}} 

\newcommand		{\bs}				{\bigskip}

\newcommand		{\mn}				{\mspace{-2mu}}
\newcommand		{\mnn}				{\mspace{-1mu}}
\newcommand		{\dsp}			{\displaystyle}

\newcommand		{\nd}				{\noindent}
\newcommand		{\ol}				{\overline}
\newcommand		{\os}			{\overset}
\newcommand		{\us}			{\underset}

\newcommand		{\ul}			{\underline}
\newcommand		{\wh}			{\widehat}

\newcommand		{\wt}			{\widetilde}

\newcommand		{\mr}			{\mathrm}
\newcommand		{\bb}			{\mathbb}
\newcommand		{\mb}			{\mathbf}

\newcommand		{\f}			{\mathfrak}
\newcommand		{\ms}			{\mathscr}

\newcommand		{\g}			{\gamma}
\newcommand		{\e}			{\epsilon}

\renewcommand	{\t}			{\theta}
\newcommand		{\s}			{\sigma}
\newcommand		{\es}			{\varsigma}
\newcommand		{\vp}			{\varphi}

\newcommand		{\G}			{\Gamma}
\newcommand		{\D}			{\Delta}

\DeclareSymbolFont{cmletters}{OT1}{cmr}{m}{n}
\DeclareMathSymbol{\Ups}{\mathalpha}{cmletters}{"7}
\renewcommand	{\Upsilon}{\Ups}

\newcommand		{\0}		{\varnothing}

\newcommand		{\F}		{\bb F}

\newcommand		{\Z}		{\bb Z}
\newcommand		{\Q}		{\bb Q}
\newcommand		{\R}		{\bb R}
\newcommand		{\C}		{\bb C}

\newcommand 		{\HH}	{\mathbb{H}}

\DeclareMathOperator{\id}			{id}

\renewcommand 		{\H}		{H^*}

\newcommand 		{\HT}	{{\H_T}}

\newcommand		{\HSN}	{H_S^N}

\newcommand		{\HKN}	{H_K^N}

\let\union\cup%
\renewcommand	{\cup}		{\mspace{-1mu}\smile\mspace{-1mu}}
\let\inter\cap%
\renewcommand	{\cap}		{\mspace{-1mu}\frown\mspace{-1mu}}
\newcommand		{\Union}		{\bigcup}

\newcommand		{\less}		{\setminus}
\newcommand		{\sub}		{\subseteq}

\newcommand		{\quot}		{\,/ \mn\mn /\,}

\renewcommand	{\-}		{^{-1}}
\renewcommand	{\o}		{\circ}
\renewcommand	{\.}		{\cdot}
\newcommand		{\x}		{\times}

\newcommand{\oplushigher}{\mathbin{\raisebox{.85pt}{$\displaystyle\oplus$}}}
\newcommand{\otimeshigher}{\mathbin{\raisebox{.85pt}{$\displaystyle\otimes$}}}
\DeclareMathOperator*{\otimesvariable}{%
			\mathchoice {\raisebox{.85pt}{$\displaystyle\otimes$}}
						{\raisebox{.85pt}{$\otimes$}}
						{\raisebox{0.7pt}{$\scriptstyle\otimes$}}
						{\raisebox{0.2pt}{$\scriptscriptstyle\otimes$}}
						}
\newcommand		{\tensor}		{\otimesvariable}
\newcommand		{\xt}[3]		{{#2}\us{#1}\otimeshigher{#3}}
\newcommand		{\direct}		{\oplushigher}
\newcommand		{\ox}			{\tensor}

\newcommand		{\+}			{\direct}
\newcommand		{\Tensor}		{\bigotimes}
\newcommand		{\Direct}		{\bigoplus}

\newcommand		{\limit}		{\varprojlim}
\newcommand		{\colim}		{\varinjlim}

\DeclareMathOperator{\diag}		{diag}

\DeclareMathOperator{\rk}			{rk }
\DeclareMathOperator{\im}		{im }

\DeclareMathOperator{\Tor}		{Tor}

\DeclareMathOperator{\chern}		{ch }

\DeclareMathOperator{\ad}		{ad }
\DeclareMathOperator{\Ad}		{Ad }

\DeclareMathOperator{\Hom}		{Hom}

\DeclareMathOperator{\Aut}		{Aut }
\DeclareMathOperator{\Inn}		{Inn }

\newcommand		{\SO}		{\mr{SO}}

\newcommand		{\U}			{\mr{U}}
\newcommand		{\SU}			{\mr{SU}}
\newcommand		{\PSU}		{\mr{PSU}}
\newcommand		{\Sp}			{\mr{Sp}}
\newcommand		{\Spin}		{\mr{Spin}}

\newcommand		{\longto} 		{\longrightarrow}
\newcommand		{\lt}			{\longto}

\newcommand		{\lmt}			{\longmapsto}

\newcommand		{\from}		{\leftarrow}
\newcommand		{\longfrom}	{\longleftarrow}

\newcommand		{\inc}		{\hookrightarrow}
\newcommand		{\xinc}		{\xhookrightarrow}
\newcommand		{\longinc}		{\xinc[]{\ \ \ \ }}

\newcommand		{\longmono}	{\xmono[]{\ \ \ \ }}

\newcommand		{\longepi}	{\xepi[]{\ \ \ \ }}

\newcommand		{\simto}		{\xrightarrow{\sim}}

\newcommand		{\longsimto}	{\os\sim\longto}
\newcommand		{\isoto}		{\longsimto}
\DeclareRobustCommand{\longsimfrom}	{\os\sim\longfrom}
\newcommand		{\isofrom}		{\longsimfrom}

\newcommand		{\ceq}			{\coloneqq}
\newcommand		{\eqc}			{\eqqcolon}

\newcommand		{\normal}			{\unlhd}

\newcommand		{\iso}				{\cong}
\newcommand		{\homeo}			{\approx}

\usepackage{todonotes}
\usepackage[russian,UKenglish]{babel}
\renewcommand\dateUKenglish{\def\today{\number\day~%
 \ifcase \month \or January\or February\or March\or April\or May\or June\or
   July\or August\or September\or October\or November\or December\fi\space
 \number\year}}
\dateUKenglish
\renewcommand{\C}{\mathbb{C}}
\renewcommand{\.}{\cdot}
\renewcommand{\d}{\delta}
\renewcommand{\b}{\beta}
\renewcommand{\H}{H^*}
\renewcommand{\f}{\mathfrak}
\let\rk\relax\DeclareMathOperator{\rk}{rk}
\usepackage{cancel}
\numberwithin{equation}{section}
\renewcommand{\U}{\mathrm{U}}
\usepackage{xcolor}
\colorlet{DarkEmerald}{Emerald!90!black}
\newcommand{\num}[1]{}
\newcommand{\rev}[1]{#1}
\newcommand{\dead}[1]{}

\theoremstyle{definition}
\newtheorem{construction}[table]{Construction}
\newextarrow{\xbigtoto}{{20}{20}{20}{20}}
{\bigRelbar\bigRelbar{\bigtwoarrowsleft\rightarrow\rightarrow}}

\newcommand{\unit}{\mathbf{1}}

\newcommand{\PSp}{\mathrm{PSp}}

\renewcommand {\^}{\wedge}

\renewcommand{\ring}{\Q}

\renewcommand{\HH}{H^*_H}

\newcommand{\KH}{K^*_H}

\newcommand{\HKSS}{Hodgkin--K\"unneth spectral sequence\xspace}
\newcommand{\KSS}{K\"unneth spectral sequence\xspace}
\newcommand{\AHSS}{Atiyah--Hirzebruch spectral sequence\xspace}
\newcommand{\AHLSSS}{Atiyah--Hirzebruch--Leray--Serre spectral sequence\xspace}
\renewcommand{\ring}{k}

\newcommand{\kWEF}{$\kk$--weakly equivariantly formal\xspace}
\newcommand{\WEFity}{weak equivariantly formality\xspace}

\newcommand{\EC}{equivariant cohomology\xspace}
\newcommand{\BEC}{Borel \EC}

\newcommand{\Tss}{\Tor^{*,*}}

\newcommand{\HHp}{H^{**}_H}

\newcommand{\thy}{ h^*}
\newcommand{\iact}{isotropy action\xspace}
\newcommand{\IDEAL}{\smash{\wh I}}
\newcommand{\etalchar}[1]{$^{#1}$}

\title{K-theory and formality}
\author{Jeffrey D. Carlson\thanks{\ %
		 	The author was partially supported by a postdoctoral fellowship 
			from the Instituto Nacional de Matem\'{a}tica Pura e Aplicada
			(IMPA)
			funded by the 
			Coordena{\c{c}}{\~a}o de Aperfei{\c{c}}oamento 
			de Pessoal de N\'{i}vel Superior (CAPES).
			}
		}
\begin{document}

\maketitle

\begin{abstract}
	We compare several notions of equivariant formality
	for complex K-theory with respect to a compact Lie group action, 
	including surjectivity of the forgetful map and
	the weak equivariant formality of Harada--Landweber,
	and find all are equivalent under standard hypotheses.
	
	As a consequence we \dead{compute} \num{8}
	\rev{present an expression for}
	the equivariant K-theory of the isotropy action of $H$ 
	on a homogeneous space $G/H$
	in all the classical cases.
	The proofs involve mainly homological algebra and 
	arguments with the \AHLSSS, but a more
	general result depends on a map of spectral sequences from
	Hodgkin's \KSS in equivariant K-theory to that in Borel cohomology
	that seems not to have been otherwise defined.

	The hypotheses for the main structure
	result are analogous to a previously announced
	characterization of cohomological 
	equivariant formality, first proved here,
	expanding on results of Shiga and Takahashi.
\end{abstract}

\section{Introduction}
The computation of the cohomology and K-theory of a homogeneous space $G/H$ 
of a compact Lie group $G$ is a venerable problem in topology.
The rational cohomology ring has been known since Cartan~\cite{cartan1950transgression}, 
and the additive structure over other coefficient rings has been recovered under reasonable hypotheses
by many authors~\cite{borelthesis,baum1968homogeneous,munkholm1974emss,gugenheimmay,husemollermoorestasheff1974,wolf1977homogeneous}, though the multiplicative structure was only confirmed \dead{very} recently~\cite{franz2021homogeneous}.
Corresponding answers in K-theory, under stronger hypotheses, came later%
~\cite{hodgkin1975kunneth,harris1968homogeneous,%
	snaith1971homogeneous,pittie1972homogeneous,minami1975symmetric};
see Remark 6.3 in this paper's prequel~\cite{carlsonfok2018} for a summary.

There are strong results 
on the structure of the equivariant K-theory ring of a compact Lie 
group action with full-rank isotropy groups, at least rationally~\cite{ademgomez2012},
but much less is known in general if the isotropy is of not of full rank.
A natural action with smaller isotropy groups to consider, then, is the so-called 
\emph{isotropy action} of $H$ on $G/H$
restricted from the transitive action,
which can be seen as a global version of the isotropy action of $H$ on the tangent space at an $H$-fixed point of a smooth $G$-manifold.
%
There is a computable description of the Borel equivariant cohomology
$H_K^*(G/H;\Q)$~\cite{carlson2016grassmannian}
for $K$ another closed, connected subgroup of $G$,
but there has been no analogous result in equivariant K-theory,
despite equivariant K-theory being a main tool in the computation of $\K(G/H)$.

In some sense the most comprehensible actions 
from an algebro-topological viewpoint are
the \emph{equivariantly formal} actions,
those for which the natural map $H_G^* X \lt H^*X$ 
from Borel equivariant to singular cohomology is surjective~\cite[Ch. XII]{borel1960seminar}%
\cite[{\SS}1.2, Thm.~14.1]{GKM1998}.
This classical
hypothesis shows up throughout
algebraic, symplectic, differential, and toric geometry,
sometimes as a tool, elsewhere as an object of study~%
\cite{BV1982,AB1984,jeffreykirwan1995,GGK,bahrisarkarsong,amannzoller2019}.
Some effort has thus gone towards determining when the
isotropy action is equivariantly formal~\cite{shiga1996equivariant,
			brion1998equivariant,
			carlson2018eqftorus,
			carlsonfok2018}.
Analogous questions have also
been explored for equivariant K-theory~\cite{matsunagaminami,haradalandweber2007},
where they amount to the even more natural question 
of whether each complex vector bundle over a $G$-space
admits a $G$-equivariant stabilization.
Fok~\cite
{fok2017formality} 
showed equivariant formality with rational coefficients
is equivalent to surjectivity of the forgetful map $\KG(X;\Q) \lt \K(X;\Q)$
from the $\Q$-localization of complex equivariant K-theory, 
and jointly with the present author~\cite{carlsonfok2018}
exploited this equivalence 
to derive a number of equivalent characterizations of \eqfity of an isotropy action.
That work did not however
compute 
$\KH(G/H\rev{;\kk})$,
for coefficient rings $\kk$ smaller than $\Q$, 
leaving that task for this sequel.

\begin{restatable}{thmx}{MAIN}\label{thm:main}
	Let $(G,H)$ be a pair of compact, connected Lie groups 
	and $k$ a subring of $\Q$ 
	such that
	$\pi_1 G$ is free abelian,
	$R(H;\kk) = RH \ox \kk$ is flat
	over the image 
	of the restriction of extended representation rings $R(G; \kk) \lt R(H ; \kk)$,
	and the kernel of this map is generated by a regular sequence.
	Then there is a $\Z/2$-graded $R(H;\kk)$-algebra isomorphism
	\[
		\KH(G/H;k) 
			\,\mn\iso\,\mn
		\xt{R(G;\kk)}{R(H;\kk)\mnn}{\mnn R(H;\kk)} \,\mnn\ox_\kk \ewP\mathrlap,
	\]
	where $\ewP \iso \im\mn\big(\K(G/H;\kk) \to \K(G;\kk)\mnn\big)$ 
	is an exterior $\kk$-algebra on $\rk_\kk \wP = \rk G - \rk H$ generators,
	$RH$ is an $RG$-module by restriction of complex representations,
	and the $R(H;\kk)$-module structure on $\K(G/H;\kk)$
	is identified with the $R(H;\kk)$-module structure on the left tensor factor.
	Similarly, there is a $\Z/2$-graded ring isomorphism
	\[
		\K(G/H;\kk) 
			\iso 
		\xt{R(G;\kk)}{\kk\mnn}{\mnn R(H;\kk)} \,\mnn\ox_\kk \ewP
	\]
	such that the forgetful map $f\: \KH(G/H;\kk) \lt \K(G/H;\kk)$
	is identified with reduction modulo the augmentation ideal $I(H;\kk)$.
\end{restatable}

As is usual with K-theory and rational cohomology,
the K-theoretic structure theorem formally resembles the cohomological one%
~\cite{carlson2016grassmannian},
but tends to hold already over $\Z$ before any primes
are inverted and requires 
a very different proof.

The hypotheses of \Cref{thm:main}
imply, by the results of the predecessor work~\cite{carlsonfok2018}, 
that the forgetful 
map $\KH(G/H;\Q) \lt \K(G/H;\Q)$
is surjective,
but the consequent 
implies the \emph{a~priori} stronger statement,
called \emph{weak equivariant formality}
by Harada--Landweber~\cite[Def.~4.1]{haradalandweber2007},
that $\K(G/H;\kk)$ 
is isomorphic to the reduction
$\KH(G/H;\kk) \ox_{R(H; \kk)} \kk$.
This turns out not to be an accident. 
Fok~\cite{fok2017formality} 
established the equivalence of
these notions over $\Q$,
but with the additional assumption $\pi_1 G$ is free abelian,
we are able, somewhat surprisingly,
to extend the equivalence to~$\Z$ coefficients.

\begin{restatable}{thmx}{WEAK}\label{thm:weak}
	Let $G$ be a compact, connected Lie group, 
	$X$ a compact $G$-space \rev{homotopy equivalent to a finite CW complex},
	and $\kk$ a subring of $\Q$.
	If $\pi_1 G$ is not free abelian,
	assume additionally the common torsion primes of
	$H^*(G;\Z)$ and $\K X$ are invertible in $k$.
	If $\KG(X;\kk)$ is finitely generated over $R(G;\kk)$
	and the forgetful map $f\: \KG(X;\kk) \lt \K(X;\kk)$ is surjective,
	then 
	the action is \kWEF.
\end{restatable}

The key tool for the proof of \Cref{thm:weak}
is an isomorphism involving two {\AHLSSS}s of fibrations 
over the base $BG$, which is of infinite type;
although this case does not satisfy the standard finiteness hypotheses, 
such sequences always strongly converge 
nevertheless (\Cref{thm:AHLSSSML}).
Our proof of \Cref{thm:main} requires the use of \Cref{thm:weak}. 
Under more general hypotheses,
\Cref{thm:condcollapse} recovers the associated graded 
ring of $\K_H(G/H;\kk)$ with respect to a certain filtration,
and the key tool for this proof is
a map of {\KSS}s:

\begin{restatable*}{theorem}{SPECMAP}\label{thm:eqkunnethchern}
Let $G$ be a compact Lie group.
The equivariant Chern character $\chern_G\: \KG \lt H_G^{**}$
induces a natural transformation of $(\Z \x \Z/2)$-graded K\"{u}nneth spectral sequences,
viewed as functors on the category of pairs $(X,Y)$ of 
compact, locally $G$-contractible $G$-spaces with finite covering dimension.
The map of $E_2$ pages is 
\eqn{
	\Tor_{RG}\blbl\mn\big(\KG X,\KG Y\big)	
	&\lt
	\Tor_{H_G^{**} \mnn}\blbl\mn\big(H_G^{**} X,H_G^{**} Y\big).\\
	\intertext{If $\pi_1 BG = \pi_0 G = 0$ then the right sequence 
		strongly converges to $H_G^{**}(X \x Y)$. 
		If $\pi_1 G$ is torsion-free as well, 
		then the left sequence strongly converges to $\KG(X \x Y)$, 
		so the $\Ei$ page is}
	\gr \KG(X \x Y)
	&\lt
	\gr H_G^{**}(X \x Y).
}
\end{restatable*}

The domain spectral sequence in this statement 
is also the main tool in the proof of \Cref{thm:main}.
This spectral sequence was defined 
by Hodgkin, but the full proof of convergence to the intended target, 
completed over a decade later, 
required the efforts of several authors and 
is somewhat circuitous,
so we present a summary in logical order in \Cref{rmk:convergence}.
We also show via \Cref{ex:SO(3)}, 
which seems not to otherwise be \rev{present} in the literature,
that the convergence of the \HKSS is not as one would hope
if $\pi_1 G$ does have torsion,
even over $\Q$ and even if $G$ is path-connected.

That the seemingly strong hypotheses to \Cref{thm:main}
are natural is demonstrated through the following corollaries,
which cover most cases where the isotropy
action of $H$ on $G/H$ is known to be equivariantly formal over $\Q$.
Moreover, by \Cref{thm:weak},
the actions are weakly equivariantly formal,
so that applying the functor $- \otimes_{RH} \Z$
recovers the ordinary K-theory $\K(G/H)$,
which was already known in these cases.

\begin{restatable*}{corollary}{TNHZ}\label{thm:TNHZ}
	Let $(G,H)$ be a pair of compact, connected Lie groups such that
	$\pi_1 G$ is free abelian.
	If $RG \lt RH$ is surjective with kernel given by a regular sequence, 
	then there are $\Z/2$-graded ring isomorphisms
	\[
	\KH(G/H) \iso RH \ox \ewP\mathrlap,\qquad\qquad\qquad
	\K(G/H) \iso \ewP\mathrlap.
	\qquad
	\]
\end{restatable*}

\begin{restatable*}{corollary}{EQUALRANK}\label{thm:equalrank}
	Let $(G,H)$ be a pair of compact, connected Lie groups such that
	$\pi_1 G$ is free abelian.
	If $H$ is of the same rank as $G$, 
	then there are $\Z/2$-graded ring isomorphisms
	\[
	\KH(G/H) \iso \xt{RG}{RH}{RH}\mathrlap,\qquad\qquad\qquad
	\K(G/H) \iso \xt{RG}{RH}{\Z}\mathrlap.
	\]
\end{restatable*}
\begin{restatable*}{corollary}{GENSYMM}\label{thm:gensymm}
	Let $(G,H)$ be a pair of compact, \dead{simply-}connected Lie groups\dead.
	\rev{such that either $G$ is simply-connected
	or the fundamental groups of both $G$ and
	\rev{the identity component of} $H$
	are free abelian.}
	If there is a continuous finite-order 
	automorphism $\s$ of $G$ 
	such that $H$ is 
	\rev{a union of components}
	\dead{the identity component}
	 of its fixed point set,	
	then there are $\Z/2$-graded ring isomorphisms
%
%
	\[
		\KH(G/H) \iso \xt{RG}{RH}{RH} \ox \ewP
		\mathrlap,\qquad\qquad\qquad
		\K(G/H) \iso 	\xt{RG}{RH}{\Z}\ox \ewP\mathrlap.
	\]
\end{restatable*}

While our results cover the best-known cases,
the regularity assumption in \Cref{thm:TNHZ} 
seems undesirable
because it is redundant in the analogous statement for rational cohomology.
One would hope it was redundant in this case also,
but it seems impossible to avoid
and is linked to the unresolved status 
of the long-standing
Abhyankar--Sathaye embedding conjecture 
in affine algebraic geometry,
discussed briefly in \Cref{rmk:Abhyankar}.

%
%
%

The cohomological analogues of our strong hypotheses are equivalent to equivariant formality~\cite{carlsonfok2018}.
This is the content of the
following extension of a theorem of Shiga and Takahashi, 
announced without proof 
in this paper's predecessor~\cite[Rmk.~3.12]{carlsonfok2018}
and finally produced in the \dead{final} \rev{concluding} \Cref{sec:ST}.

\setcounter{thmx}{2}
\begin{restatable}{thmx}{STPPP}\label{thm:ST+++}
Let \GK be a pair of compact, connected Lie groups.
\TFAE: 
\begin{enumerate}
\item[\textup{(1)}] The group $K$ acts equivariantly formally on $G/K$.
\item[\textup{(2)}] The space $G/K$ is formal
		and the map $H^*(BG;\Q) \lt H^*(BK;\Q)^{N_G(K)}$ is a surjection.
\item[\textup{(3)}] The finite group $\pi_0 N_G(K)$ acts on the space of 
		indecomposables $QH_K^*$ as a reflection group and
		the map $H^*(BG;\Q) \lt H^*(BK;\Q)^{N_G(K)}$ is a surjection.
\item[\textup{(4)}] The finite group $\pi_0 N_G(S)$ acts on the 
		tangent space $\fs$ to a maximal torus $S$ of $K$ as a reflection group and 
		the map $H^*(BG;\Q) \lt H^*(BS;\Q)^{N_G(S)}$ is a surjection.
\end{enumerate}
\end{restatable}

\medskip

The structure of the paper is as follows.
In \Cref{sec:torsion}\rev{,} \dead{contains the proof of} \rev{we prove} \Cref{thm:weak}
and develop\dead{s} a number of other equivalences 
of equivariant formality--like notions.
In \Cref{sec:main}, 
we prove \Cref{thm:main}, \Cref{thm:condcollapse}, and their corollaries.
In \Cref{sec:mapkunneth} we set up 
the spectral sequences of \Cref{thm:eqkunnethchern}
and the map therebetween.
Finally, in \Cref{sec:ST} we prove the expanded Shiga--Takahashi \Cref{thm:ST+++}. 


\medskip

\nd 
\emph{Acknowledgments.} 
The author is grateful to Omar Antol\'{i}n Camarena and Larry Smith 
for many helpful and enjoyable conversations
and for the generosity of the  
National Center for Theoretical Sciences in Taipei,
where a portion of this work was conducted.

\section{Equivariant formality and torsion}\label{sec:torsion}

As mentioned in the introduction, 
equivariant formality for rationalized K-theory and for 
Borel cohomology are the same phenomenon~\cite{fok2017formality} 
and there are by now many equivalent formulations~\cite{carlsonfok2018}. 
In this section, 
we add another and explore to what extent the equivalence 
with weak equivariant formality, to be discussed, \dead{can} survive\rev{s} with more general coefficients.

\begin{notation}
	In all that follows, $\defm G$ will be a compact, connected Lie group.
\end{notation}
	
\begin{definition}
	
	The \defd{equivariant K-theory}
	$\defm{K^0_G X}$ of a \rev{compact} $G$-space $X$ \num{Def. 2.2}
	is the additive Grothendieck group of the semiring of 
	$G$-equivariant vector bundles $V \to X$, 
	bundles whose total space admits a $G$-action such that the 
	projection is equivariant.
	%
	The unique $G$-map $X \to {*}$ induces a map of 
	groups $K_G^0({*}) \longrightarrow K_G^0(X)$ 
	whose cokernel is the \emph{\textbf{reduced equivariant K-theory}}
	$\widetilde{K}^0_G X$.
	We set $\widetilde{K}^{-n}_G X = K^0_G \Sigma^n X$ for $n \geq 0$
	and $\Sigma$ the reduced suspension,
	and $K^{-n}_G X = \widetilde{K}^{-n}_G(X_+)$
	for $X_+$ the union of $X$ and a disjoint, $G$-fixed basepoint $+$.
	As with
	nonequivariant K-theory, Bott periodicity holds,
	and we can form a $\Z/2$-graded ring $\defm{\KG X} \ceq K^0_G X \+ K^1_G X$.
	For each characteristic-zero coefficient ring $k$,
	the functor
	given on finite $G$--CW-complexes by
	$X \lmt \KG X \ox k$
	extends uniquely to an equivariant cohomology theory		
	$\defm{\KG(-;k)}$.\footnote{\ 
		\rev{It is uncommon to see $\KG X$ considered 
			for non-compact spaces, 
			and in this paper we continue the tradition of not discussing it,
			but at the urging of the referee,
			we provide some variety by including an explicit definition.
			In any event we set 
			$K^{-1}_G(X) = \wt K^0_G(\susp X)$.
			For locally compact Hausdorff spaces $X$,
			one obtains a workable compactly-supported cohomology theory
			by considering the one-compactification $X^+$,
			with compactifying point $+$,
			and setting $\wt K_G(X) \ceq \wt K^*_G(X^+,+)$.
			The correct homotopical definition~\cite[{\S}XIV.4]{alaska}
			generalizes the nonequivariant definition 
			$K^0(X) = [X,B\U\x \Z]$ 
			and
			$K^{-1}(X) = [X,\U]$.
			Consider the direct sum $\defm{\mathbf U} \ceq \Direct W^\infty$ 
			of $\aleph_0$-many
			copies of each irreducible $G$-representation $W$,
			equipped with a $G$-invariant Hermitian inner product.
			For each $n$ and finite-dimensional $G$--vector subspace $V < \mb U$,
			the Grassmannian $G_n(V \+ \mb U)$ of 
			complex $n$-planes in $V \+ \mb U$ inherits a $G$-action,
			and inclusions $V \leq V' = V^\perp \+ V$ induce
			maps $G_n(V\+\mb U) \lt G_{n+|V^\perp|}(V'\+ \mb U)$
			via $L \lmt V^\perp \+ L$.
			We set $B\U_{2n}^G(V) \ceq \coprod_n G_{2n+|V|}(V+\mb U)$,
			so that inclusions $V \leq V'$ induce
			basepoint-preserving $G$-maps $B\U^G(V) \lt B\U^G(V')$,
			and 
			let $B\U^G \ceq \colim_V B\U^G(V)$,
			and then $\wt K^0_G(X) = [X,B\U^G]_G$,
			where $[-,-]_G$ denotes basepointed $G$-homotopy classes
			of basepointed $G$-maps.
		}
	}
	\rev{For $G = 1$ this recovers the classical homotopical definition of }
\num{Def. 2.2}
	\rev{nonequivariant K-theory: $\wt K^0(X) \ceq [X,B\U]$ and $K^1(X) = [X,\U]$.}
	An equivariant $G$-bundle over a single point $*$ is just a representation, 
	so $\KG(*)$ is the complex \defd{representation ring} $\defm{RG}$ of $G$.
	This ring carries a canonical augmentation 
	$\defm{\e}\: RG \lt \Z$
	defined on a representation $\rho\: G \lt {\Aut}_{\C} V$ by $\rho \mapsto \dim V$,
	whose kernel is written $\defm{IG}$.
	We also write $\defm{R(G;\kk)} = RG \ox \kk$, and $\defm{I(G;\kk)} \iso IG \ox k$.
	The unique $G$-map $X \to {*}$
	induces~\cite[\SS4.5]{atiyahhirzebruch} 
	a map $R(G;k) \lt \KG(X;k)$
	which makes $\KG(X;k)$ an algebra over $R(G;k)$
	just as $H_G^{**} X$ is an algebra over $H_G^{**}$.
	Forgetting the $G$-structure on a bundle induces a natural \defd{forgetful map}%
	~\cite{matsunagaminami}
	\[ 
	\defm{f}\: \KG X \lt \K X.
	\]
	Surjectivity of $f \ox \id_\Q\: \KG(X;\Q) \lt \K(X;\Q)$
	is called (\emph{K-theoretic}) \defd{\eqfity},
	and \rev{is} equivalent~\cite[Thm.~1.3]{fok2017formality}\cite[Thm.~5.6]{carlsonfok2018} \num{143}
	to the map $X \longinc X_G$
	inducing a surjection in cohomology.
\end{definition}
\begin{definition}
	We denote a universal principal $G$-bundle by $EG \to BG$ 
and the homotopy orbit space of a continuous action of $G$ on a space $X$ 
by $\smash{\defm{X_G} \,\ceq\, EG \ox_G X 
	\ \mnn = \ \mnn 
	\quotientmedddd{EG \x X\mspace{3.25mu} }{\mspace{3.25mu} (eg,x) \sim (e,gx)}}$.
For us, \emph{\textbf{cohomology will take rational coefficients}} by default,
and we will more frequently need to view it as the  
\emph{product} $\defm{\Hp} \ceq \prod H^n$ rather than the
  coproduct $\defm{H^*}\ceq \Direct H^n$ of its graded components.
By definition, the \defd{Borel equivariant cohomology} of a continuous action
of $G$ on a space $X$
is the singular cohomology $\defm{H_G^{**} X} \ceq \Hp X_G$
of the homotopy orbit space.
It will be convenient at times to also write $\defm{H_G^{**}} \ceq H_G^{**}(*) = \Hp BG$
for its coefficient ring
and likewise $\defm H_G^* = H_G^*(*) = H^*BG$.
\end{definition}

\begin{discussionsilent}
We make essential use of the 
natural transformation from equivariant K-theory to Borel cohomology
given by the \defd{equivariant Chern character} 
$\defm{\chern_G}\:\KG \lt H_G^{**}$,
defined on the class of an equivariant vector bundle $V \to X$ by
\eqn{
	\KG X \lt 	&\,\ \, \K X_G \, \, \lt \Hp X_G,\\
	V \lmt 		&\, \xt G {EG} V \lmt \chern(\xt G {EG} V).
}
The first map is completion, by the theorem of Atiyah and Segal~%
\cite[Prop 4.2]{atiyahsegalcompletion}, 
and the second the ordinary Chern character,
which is a 
natural
$\Z/2$-graded ring isomorphism $\K(X;\Q) \lt \Hp X$.
Particularly, we find~\cite[Thm.~5.3]{carlsonfok2018}
that $\chern_G$ induces isomorphisms 
\quation{\label{eq:extension}
	\KG(X;\Q)\compl \isofrom \xt{RG}{\KG X}{H_G^{**}} \isoto H_G^{**} X
}
natural in finite $G$--CW complexes $X$.
\end{discussionsilent}

The Chern character is frequently injective.
The following standard observation will play a crucial role.

\begin{lemma}\label{thm:extflat}
	If $G$ is a compact, connected Lie group, 
	the composition of the following ring extensions is a flat injection:
	\[
	RG 
	\longmono 
	\wR G  
	\longmono 
	\wR G  \ox \Q 
	\longmono
	R(G;\Q)\compl 
	\isoto\,
	{H_G^{**}}.
	\]
\end{lemma}
\bpf
The penultimate map is a completion and the last is \eqref{eq:extension} 
for $X = {*}$.
For flatness, $\Q$ is flat over $\Z$,
so $\wR G \ox \Q$ is flat over $\wR G $,
and $RG \lt \wh RG$ and $\wh RG \ox \Q \lt (\wh RG \ox \Q)\compl$ 
are flat extensions since $RG$ is Noetherian%
~\cite[Cor.~3.3]{segal1968representation}\cite[Prop.~10.14]{atiyahmacdonald}
and hence $\wh RG$ is Noetherian~\cite[Thm.~10.\rev{26}\dead{14}]{atiyahmacdonald}.
For injectivity, factor the sequence as
$RG \to R(G;\Q) \to R(G;\Q)\compl$ instead.
The first map is injective because $RG$ 
is free abelian on the irreducible representations; 
the second is the restriction to
Weyl group invariants of the completion $R(T;\Q) \lt R(T;\Q)\compl$,
which is the embedding
$\Q[t_1^{\pm_1}, \ldots, t_n^{\pm 1}] \longinc \Q[[t_1 - 1,\ldots,t_n-1]]$%
~\cite[Thm.~4.4, Prop.~4.3]{atiyahhirzebruch}.
%
%
\epf
\brmk\label{rmk:completion-weak}
The sequence $RG \to \wh{R}G \overset\sim\to K^*BG \to H_G^{**}$
was already discussed by Atiyah and Hirzebruch  \cite[\SS4.7]{atiyahhirzebruch}.
The map $K^*BG \otimes \mathbb Q \lt H_G^{**}$ is
not yet an isomorphism,
for example, when $G = S^1$\rev{,} because $\mathbb Q[[u]]$ 
is much larger than the subring $\mathbb Z[[u]] \otimes \mathbb Q$
of power series with coefficients of bounded denominator. 
The isomorphism of \eqref{eq:extension} yields less information about
K-theory than might first appear; 
not only does it annihilate torsion, 
but completion is so blunt an instrument that, for example,
the completion with respect to the augmentation ideal of
each of the three pairwise nonisomorphic rings 
$R\SU(3) \iso \Z[x,y]$ 
and $R\PSU(3) \iso \Z[a,b,c]/(a^3-bc)$
and $R\U(1)^2 \iso \Z[t^{\pm 1},u^{\pm 1}]$
is a power series ring in two indeterminates.
\ermk

\begin{discussionsilent}
	The fact that equivariant formality is equivalent to the collapse of the
	\SSS of the Borel fibration $X \to X_G \to BG$ 
	and to a reduction isomorphism $H^*X \iso H_G^* X \ox_{H_G^*} \Q$
	suggests a variant notion of equivariant formality.
	Note that the augmentation ideal $IG$
	is generated by differences $\rho - \e(\rho)$,
	so that the image of $IG$ under the structure map $RG \lt \KG X$
	lies in the kernel of the forgetful map
	$f\: \KG X \lt \K X$.
	Thus $f$  factors through the quotient
	$
	\defm{{\KG X} \quot {RG}} \,=\, 
	{\KG X}\ox_{RG}{\Z},
	$
	where $\e$ provides the $RG$-module structure on $\Z$,
	and we have a diagram
	\quation{\label{eq:HLdiagram}
		\begin{aligned}
			\xymatrix@C=0.5em@R=1em{
				R\dead(G\dead) \ar[rr]& \ & \KG X \ar[rr]^f \ar@{->>}[rd] && \K X \\
				&&& \qquotientmedd{\KG X}{ R\dead(G\dead) } \ar[ur]_(.6){\defm{\bar f}} &
			}
		\end{aligned}
	}
	analogous to the cohomology $H_G^* \to H_G^* X \to H^*X$ of the Borel fibration.
	Harada and Landweber
	observe 
	that $f$ is surjective if and only if $\bar f$ is%
	~\cite[Prop.~4.2]{haradalandweber2007},
	and moreover that if $\bar f$ is an {isomorphism},
	then $f$ is surjective with kernel $IG \. \KG X$.
	They thus set the following definition.
\end{discussionsilent}

\begin{definition}[{\cite[Def.~4.1]{haradalandweber2007}}]
	\label{def:weakKEF} 
	Let $\kk$ be a ring.
	A $G$-action on a space $X$ is \defd{$\kk$--\WEF} 
	if the forgetful map induces an isomorphism
	\[
	\qquotientmed{\KG(X;\kk)} {R(G;\kk)} 
	\,\lt\,
	\K(X;\kk).
	\]
	We simply say the action is \defd{\WEF}
	in the case $\kk = \Z$.
\end{definition}

\begin{discussionsilent}
	We have these evident implications for $\kk$ a subring of $\Q$:\num{189}
	\[
	\xymatrix{
		\big[{\KG(X;\kk)} \longepi {\K(X;\kk)}\big] \quad
		\ar@{=>}[r]  &
		\quad\big[{\K(X_G;\kk)} \longepi {\K(X;\kk)}\big]
		\ar@{=>}[d]
		\\
		\Big[\qquotientmed{\KG(X;\kk)} {R(G;\kk)} \iso \K(X;\kk)\Big]\quad
		\ar@{=>}[u] 
		&
		\quad\big[{H_G^*(X;\Q)} \longepi {\dead{\K}\rev{H^*}(X;\Q)}\big] 
		\mathrlap.
	}
	\]
	Fok has shown that for $\kk = \Q$,
	the conditions are in fact all equivalent~\cite{fok2017formality}\cite[Thm.~5.6]{carlsonfok2018},
	and we examine to what extent these
	implications are reversible 
	for general subrings $\kk \leq \Q$.
	Subject to torsion hypotheses,
	we find the converses do hold:\num{192}
	\[
	\xymatrix{
		\big[{\KG(X;\kk)} \longepi {\K(X;\kk)}\big] \quad
		\ar@{=>}[d]_{\text{{\vphantom{L4}}Thm. \ref{thm:weak}\ }}
		^{\text{\ {\vphantom{L4}}Thm. \ref{thm:weak-2}}}
		\ar@{<=>}[r]^{\text{Prop. \ref{thm:KK-to-K}}{\vphantom{TW}}}&
		\quad\big[{\K(X_G;\kk)} \longepi {\K(X;\kk)}\big]
		\\
		\Big[\qquotientmed{\KG(X;\kk)} {R(G;\kk)} \iso \K(X;\kk)\Big]	\mathrlap.\quad
		&
		\quad\big[{H_G^*(X;\kk)} \longepi {\rev{\H}\dead{\K}(X;\kk)}\big]
		\ar@{=>}[u]_{\text{\Cref{thm:H-to-KK}}} 
	}
	\]
	Most of the implications are routine,
	but \Cref{thm:weak} is somewhat unexpected,
	and the proof is correspondingly more involved.
	
	We first deal with the simpler implications in order,
	taking as given, for the moment, that the 
	{\AHLSSS}s 
	\eqn{
		H^*(BG;\kk) &\implies \K(BG;\kk),\\
		H^*(X_G;\kk) &\implies \K(X_G;\kk),\\
		\H\big(BG;\K(X;\kk)\mnn\big) &\implies \K(X_G;\kk)
	}
	converge, which will follow from \Cref{thm:AHLSSSML}.
\end{discussionsilent}

\begin{lemma}\label{thm:H-to-KK}
	Let $G$ be a compact, connected Lie group, 
	$\kk$ a subring of $\Q$,
	and $X$ a $G$-space such that 
	such that $H^*(X;\kk)$ and $H^*(BG;\kk)$ are 
	free $\kk$-modules
	with $H^*(X;\kk)$ of finite rank.
	If the cohomological fiber restriction
	$g\: H_G^*(X;\kk) \lt H^*(X;\kk)$
	of the Borel fibration is surjective,
	then so is the K-theoretic fiber restriction
	$\smash{\wh f}\: \K(X_G;\kk) \lt \K(X;\kk) \vphantom{x^{X^t}}$.
\end{lemma}
\begin{proof}
	The surjectivity of $g$ implies the \SSS
	of $X \to X_G \to BG$ collapses,
	implying $H^*(X_G;\kk)$ too is free over $\kk$
	\rev{since $\kk$ is a localization of $\Z$ and hence \num{Lemma 2.12.}
	there is no extension problem for free $\kk$-modules}.
	Note that we may take CW-approximations of $X_G$ and $X$
	and homotope $X \lt X_G$ to be cellular, 
	showing the filtration on K-theory 
	implicated in the \AHSS is preserved.
	The freeness assumptions guarantee
	the {\AHSS}s $H^*(X_G;\kk) \implies \K(X_G;\kk)$
	and $H^*(X;\kk) \implies \K(X;\kk)$ collapse at $E_2$,
	\rev{as these spectral sequences are obtained  \num{Lemma 2.12.}
	from the spectral sequence over $\Z$ by tensoring with $\kk$,
	and all images of differentials in this spectral sequence are torsion;}
	\dead{so that} \rev{thus} $g$ is the associated
	graded map of $\wh f$.
	Since the filtration on $\K(X;\kk)$ is finite, 
	surjectivity of $g$ also implies surjectivity of $\wh f$.
\end{proof}

\begin{proposition}\label{thm:KK-to-K}
	Let $G$ be a compact, connected Lie group,
	$\kk$ a subring of $\Q$, 
	and $X$ a compact $G$-space  such that $\KG(X;\kk)$ is a finite $RG$-module.
	Then 
	\benum[label={({\alph*})}]
	\item\label{thm:KK-K-surj}
	$f\: \KG(X;\kk) \lt \K(X;\kk)$ is surjective
	if and only if 
	$\smash{\wh f}\:\K(X_G;\kk) \lt \K(X;\kk)$ is and
	\item\label{thm:KK-K-inj}
	$\bar f\:  \xt {\!\!R(G;\kk)\!} {\KG(X;\kk)} \kk  \lt \K(X;\kk)$ is injective
	if and only if
	$\bar{\wh f}\:\xt{\!\!\K(BG;\kk)\!}{\K(X_G;\kk)} \kk  \lt \K(X;\kk)$ is.
	\eenum
\end{proposition}
\begin{proof}
\dead{	The map forgetting the action on a
	$G$--vector bundle $V \to X$ 
	factors through the extension map~%
	$\defm i\: {\K_G X} \lt \KG (EG \x X) \iso \K(X_G)$ 
	taking a $G$--vector bundle 
	$V \to X$ to the 
	$G$--vector bundle~%
	$EG \x V \to EG \x X$,
	since the restriction of the latter bundle to the fiber over a basepoint
	${*} \in EG$ is
	$V \to X$ again.
	Thus $f$ factors as $\smash{\wh f} \o i$.
	The fiber restriction $\wh f$ annihilates the
	ideal generated by the image of~$\defm{\wh I} =  
	\K(BG,*;\kk) \to \K(X_G;\kk)$
	by definition, so it}
\rev{  
	The map 
	$\defm i\: {\K_G X} \lt \KG (EG \x X) \iso \K(X_G)$ 
	takes a $G$--vector bundle 
	$V \to X$ to the 
	$G$--vector bundle~%
	$EG \x V \to EG \x X$.
	The restriction $\wh f$ to the fiber over a basepoint
	${*} \in EG$ takes the second bundle to
	$V \to X$ again.
	Thus 
	the map $f$ 
	forgetting the action on a
	$G$--vector bundle $V \to X$ 
	factors as $\smash{\wh f} \o i$.
	The 
	ideal generated by the image of~$\defm{\wh I} =  
	\K(BG,*;\kk) \to \K(X_G;\kk)$
	lies in the kernel of $\wh f$
	by definition, so $\wh f$ 
      }
	factors through
	$\K(X_G;\kk) / \wh I \. \K(X_G;\kk)$.
	Now, $i$ 
	may be identified with completion 
	of $M = \KG(X;\kk)$ at $\defm I = I(G;\kk)$
	by the theorem of Atiyah and Segal~\cite[Cor.~2.2]{atiyahsegalcompletion},
	and $\dead{\wh I \. \K(X_G;\kk) =} \wh I \wh M$
	is \rev{also} the completion of~$I \. \KG(X;\kk) = IM$
	by standard commutative algebra.%
	\footnote{\ 
		Explicitly, 
		setting the abbrevation $\defm R = R(G;\kk)$
		as well,
		one has a string of standard isomorphisms~\cite[Props.~10.13,\,15]{atiyahmacdonald}
		\[
		\IDEAL \. \wh M
		\iso 
		\IDEAL		\otimeshigher_{\wh R}	{\wh M}
		\iso
		{(\xt{R}I{\wh R})}
		\otimeshigher_{\wh R}
		{(\xt{R}{\wh R}M)}
		\iso
		\xt{R}{\xt{R}{I}{\wh R}}{M}
		\iso
		\xt{R}{\wh R}{(I \ox_R M)}
		\iso 
		\xt{R}{\wh R}{IM}
		\iso 
		\wh{IM}
		\mathrlap.
		\]
	}
	As $R(G;\kk)$ is Noetherian and the modules in question
	are finitely generated,
	one has $\wh{M/IM} \iso \wh M / \wh{IM}$
	on completing the short exact sequence~$IM \to M \to M/IM$;
	but $M/IM$ is already complete with respect to $I$
	since $I(M/IM) = 0$,
	so~$M/IM \iso \wh M / \wh{I}\wh{M}$.
	Thus we may identify the maps ${\bar{\wh f}}$ and $\bar f$\dead{,}\rev{.
	\ref{thm:KK-K-inj} then follows immediately, 
	and 
	\ref{thm:KK-K-surj} 
	follows as well because $f$ is surjective if and only if $\bar f$ is,
	and $\wh f$ is surjective if and only if ${\bar{\wh f}}$ is.
	}
	\dead{and
	\ref{thm:KK-K-surj} (respectively, \ref{thm:KK-K-inj})
	follows since each map is 
	then surjective (resp., injective) precisely when the other is.}
\end{proof}

Weak \eqfity naturally breaks into surjectivity and injectivity clauses,
and \rev{in the free case} there is a result for injectivity too.

\begin{proposition}[Snaith {\cite[Thm.5.4]{snaith1971massey}}]%
	\label{thm:Snaith-inj}
	Let $G$ be a compact, connected Lie group,
	$X$ a compact, locally contractible $G$-space of finite covering dimension,
	and $\kk$ a subring of $\Q$
	such that
	$H^*(X_G; \kk)$ and $H^*(BG;\kk)$ are free over $\kk$ and
	the map
	$H^*(X_G; \kk) \quot H^*(BG;\kk) \lt H^*(X; \kk)$ 
	is injective.
	Then $\KG(X;\kk)\quot R(G;\kk) \lt \K(X;\kk)$ 
	is injective.
\end{proposition}
Snaith really shows this for $\kk = \Z$ but the extension
follows by tensoring his proof with $\kk$. 
\begin{corollary}\label{thm:Jeff-inj}
	Let $G$ be a compact, connected Lie group,
	$X$ a compact, locally contractible $G$-space of finite covering dimension,
	and $\kk$ a subring of $\Q$
	such that
	$H^*(X; \kk)$ 
	and $H^*(BG;\kk)$ are free over $\kk$ and
	$H^*(X_G; \kk)  \lt H^*(X; \kk)$ 
	is surjective.
	Then $\KG(X;\kk)\quot R(G;\kk) \lt \K(X;\kk)$ 
	is injective.
\end{corollary}
\begin{proof}
	The hypotheses imply
	the cohomological \SSS of $X \to X_G \to BG$
	collapses at $E_2 \iso H^*(BG;\kk) \ox_\kk H^*(X;\kk) = \Ei$
	so $H^*(X_G;\kk)$ is also free and the kernel of the edge map 
	$\Ei \longepi H^*(X;\kk)$
	is precisely the ideal generated by $H^{\geq 1}(BG;\kk)$.
	Then Snaith's \Cref{thm:Snaith-inj} applies.
\end{proof}
With the hypotheses of \Cref{thm:Jeff-inj}
in place, then, one only needs that 
$\KG(X;\kk) \lt \K(X;\kk)$ is surjective
as well to conclude $\kk$--\WEFity:
\begin{proposition}\label{thm:weak-2}
	Let $G$ be a compact, connected Lie group,
	$X$ a compact, locally contractible $G$-space of finite covering dimension,
	and $\kk$ a subring of $\Q$
	such that
	$H^*(X; \kk)$
	and 
	$H^*(BG;\kk)$ are free over $\kk$ and
	$\K_{\rev{G}}(X_{\dead{G}}; \kk)  \lt \K(X; \kk)$ \num{Prop. 2.16}
	is surjective.
	Then $\KG(X;\kk)\quot R(G;\kk) \lt \K(X;\kk)$ 
	is an isomorphism.
\end{proposition}
\begin{proof}
	Since $\kk$--\WEFity is the map 
	$\KG(X;\kk) \quot R(G;\kk) \lt \K(X;\kk)$
	being both injective and surjective
	and we hypothesize surjectivity,
	we will be done if we can demonstrate that
	$H^*(X_G; \kk)  \lt H^*(X; \kk)$ 
	is also surjective, 
	for then \Cref{thm:Jeff-inj} will apply.
	But $\K(X_G;\kk) \lt \K(X;\kk)$
	is surjective by \Cref{thm:KK-to-K},
	and localizing and completing,
	the map $\Hp(X_G;\Q) \lt H^*(X;\Q)$ is also surjective,
	or equivalently, 
	$H^*(X_G;\Q) \lt H^*(X;\Q)$ is.
	Thus the \SSS of $X \to X_G \to BG$ collapses over $\Q$.
	By the torsion hypotheses,
	the \SSS over $\kk$ 
	injects into the \SSS over $\Q$,
	which is the localization of the former,
	so if the former supported any nonzero differential,
	so also would the latter. 
	Thus the spectral sequence collapses
	and $H^*(X_G;\kk) \lt H^*(X;\kk)$ is surjective
	as desired.
\end{proof}

Thus, with sufficiently brutal hypotheses on torsion,
$\kk$--\WEFity reduces
to ordinary K-theoretic or even cohomological \eqfity.
These hypotheses are a bit unnatural, 
because, for instance, $RG$ is always torsion-free
but $H^*(BG;\Z)$ typically is not,
and one does not want to assume anything of $X$ 
beyond a finiteness condition if one can avoid it.
Fortunately, one can weaken these hypotheses substantially,
to nothing if $\pi_1 G$ is free abelian.
For this we will invoke a reduction to the maximal torus.


\begin{lemma}[Harada--Landweber~{\cite[Lem.~4.4]{haradalandweber2007}}]\label{thm:HLGT}
	Let a compact, connected Lie group $G$ with maximal torus $T$ 
	act on a finite CW complex $X$, and
	suppose $\pi_1 G$ is torsion-free.
	Then the $G$-action on $X$
	is \WEF
	if and only if the $T$-action is.
\end{lemma}

\bthm[Matsunaga--Minami~{\cite[Thm.~5]{matsunagaminami}}]
Given a compact, connected Lie group $G$ with maximal torus $T$,
there exists a positive integer $\defm m$
such that for any compact $G$-space $X$ such that $\KG X$ is a finite $RG$-module,
if $\KT X \lt \K X$ is surjective, then the image of $\KG X \lt \K X$ 
contains $m \. \K X$.
If $\pi_1 G$ is torsion-free,  then $m = 1$.
\ethm

\begin{discussionsilent}
	In the proof of \Cref{thm:weak}
	(and actually already in \Cref{thm:H-to-KK}), 
	we will use a convergence result regarding the \AHLSSS.
	This is a (cohomological, for us) right half-plane
	spectral sequence corresponding to 
	the Cartan--Eilenberg system
	(\emph{extended} in the language of Helle--Rognes~\cite[Def.~6.1]{hellerognes2018boardman})
	$H(p,q) = \E(X_{q-1},X_{p-1})$ (\cite[Ex.~2, p.~335]{cartaneilenberg}; 
	see Atiyah--Hirzebruch~\cite[Pr.,~\SS2.1]{atiyahhirzebruch}
	or reason directly for the index shift)
	for a fibration
	$\smash{F \to X \os{\defm\xi}\to }B$ be a fibration over a CW-complex (or CW-spectrum)
	with skeletal filtration $(\defm{B_p})$
	and induced filtration
	$\defm{X_p} = \xi\-B_p$.
	This spectral sequence has $E_2 = \H\big(B; \E X\big)$,
	where cohomology is considered as a $\Z[\pi_1 B]$-module;
	for more a general filtration, one 
	still has $E_1^{p,q} = E^q(X_{p},X_{p-1})$
	and the
	associated exact couple (\emph{left couple} per Helle--Rognes~
	\cite[Def.~7.1]{hellerognes2018boardman})
	has data
	\quation{\label{eq:space-filtration-SS}
		A_p = \E(X,X_{p-1}),\qquad
		\colim A_p =  \E X,\qquad
		F_p \E X = \im\!\big(  \E(X,X_{p-1}) \lt \E X\big)\mathrlap.
	}
\end{discussionsilent}

\bprop\label{thm:filt-conv}
Let $X$ be a space or connective spectrum, $(X_p)_{p \in \Z}$ 
an increasing filtration with $X_{-1} = \0$ and $X = \Union X_p$,
and $\E$ a cohomology theory.
Then 
the filtration of \eqref{eq:space-filtration-SS}
is exhaustive and complete,
and 
the associated spectral sequence 
beginning with $E_1^{p,q} = E^{p+q}(X_p,X_{p-1})$
is strongly convergent if additionally
\dead{if}\num{Prop. 2.21} each $E_1^{p,q}$ is an Artinian module over 
$\wt E^0(S^0)$ or if  
$\limit^1_p \E(X_p,X_{p_0})$ vanishes for each $p_0$.
\eprop

\bpf
Exhaustion is immediate because
$F_p \E X = \E X$ for $p < 0$.
Since $\limit A_p$ and $\limit^1 A_p$ 
are zero~\cite[Thm.~4.3(b)]{boardman1999conditionally}, 
the associated exact couple is conditionally convergent
to the colimit $\E X$
in the sense of Boardman~\cite[Def.~5.10]{boardman1999conditionally},
and hence, as a spectral sequence with entering differentials,
converges strongly to $\E X$
if we can show that $\limit^1 Z_r = 0$~\cite[Thm.~7.1 \emph{et seq.}]{boardman1999conditionally}.
Boardman observes \rev{one}
can show this separately for each $Z^{p,q}_r$,
but it is equivalent to consider a tail of this
sequence in which $B^{p,q}_r = 0$\dead,\rev.
Given our Artinian assumption, 
the decreasing sequence $Z^{p,q}_r$ will stabilize
by some finite $r$, so that $\dsp\limit_r{\!}^1 Z^{p,q}_r = 0$.

For the other hypothesis,
write $j_p\: \E(X,X_{p-1}) \lt \E X$,
so that $F_p = \im j_p$,
and consider the six-term
$\smash{\limit^\bl}$ exact sequence
associated to the 
short exact sequence
$\smash{0 \to \ker j_p \to A_p \to F_p \to 0}$.
As $\smash{\limit^1 A_p}$ vanishes, 
so too does $\smash{\limit^1 F_p}$,
so the filtration is complete.
The same exact sequence
shows $\smash{\limit F_p}$ will vanish---%
and hence, by definition, $F_p$ is Hausdorff---%
if and only if $\smash{\limit^1} \ker j_p$ vanishes.
But
writing $\d_p\: E^{*-1} X_{p-1} \lt \E(X,X_{p-1})$ for the connecting
map in the long exact sequence of the pair, 
one has a short exact sequence
$0 \to \ker \d_p \to \E X_{p-1} \to \ker j_p \to 0$,
so the $\smash{\limit^\bl}$ exact sequence
shows that $\smash{\limit^1} \ker j_p$ vanishes
so long as $\smash{\limit^1} \E X_p$ does.
Thus we will have strong convergence 
so long as we have weak convergence.
We indeed do by 
the exact sequence~\cite[p. 615--6]{whiteheadelements}
\[
0 
\,\to\, 
{\gr}{}_p \, \frac{E^{p+q} X} {{\underset r \limit}^1 E^{p+q-1} X_r} 
\,\longto\, 
E^{p,q}_\infty
\,\longto\,
\ker\mn\Big(	{\limit_r}^1 E^{p+q}(X_{p+r},X_p) 
\mn\to\,\mn
{\limit_r}^1 E^{p+q}(X_{p+r},X_{p-1})
\mnn\Big)
\,\to\,
0
\]
because we have assumed the $\limit^1$ terms vanish.
\epf

\brmk\label{rmk:ML}
To achieve the vanishing of $\limit^1$ terms 
in \Cref{thm:filt-conv},
it is sufficient to assume for each $n$
that the sequence $(E^n X_r)$ satisfies the Mitt\"ag--Leffler
criterion,
for
considering the natural maps between 
the long exact sequences of pairs $(X_{p+r},X_p)$,
a brief diagram chase factoring through the four lemma
shows that if we fix $n$, $p$, and $r_0$
and let $r$ increase,
the images of $E^n(X_{p+r},X_p) \lt E^n(X_{p+r_0},X_p)$
stabilize as soon as those of 
$E^n X_{p+r} \lt E^n X_{p+r_0}$ do.
\ermk

\bcor\label{thm:AHLSSSML}
Let $F \to X \os\xi\to B$ be a fibration over a CW-complex (or \rev{connective}
CW-spectrum)
with skeletal filtration $(B_p)$,
let $X_p = \xi\-B_p$,
and suppose for each $n$ that 
$(E^n X_p)_p$ 
is Artinian over $\wt E^0(S^0)$
or satisfies the Mitt\"ag--Leffler criterion.
Then the \AHLSSS starting with $E_2^{p,q} = H^p(B;E^q F)$ 
strongly converges to $\E X$.
\ecor

\begin{discussionsilent}
	This is a little more latitude than the more 
	commonly encountered 
	hypotheses that $B$ be finite-dimensional or that
	$E^{\leq q} F = 0$ for some $q$
	and saves us the additional step of running the spectral
	sequence separately on the skeleta $B_p$
	and collating the results through a limit argument.
\end{discussionsilent}

\WEAK*
\rev{
The plan of the proof is to use \Cref{thm:KK-to-K}
to replace the claimed isomorphism 
\[\KG(X;\kk)/I(G;\kk)\.\KG(X;\kk) \lt \K(X;\kk)\]
with the equivalent claimed isomorphism 
$\K(X_G;\kk)/\wt K^*(BG;\kk)\.\K(X_G;\kk) \lt \K(X;\kk)$,
which is susceptible to an \AHLSSS argument.
}

\rev{	
At this point the two different hypotheses lead to  
divergent middle sections of the proof.
In both cases we need to show 
that the $E_2$ page of the \AHLSSS $\big(E_r(\xi),d_r\big)$ 
of $X \to X_G \to BG$
is of a particularly simple form, 
which arises from the nonexistence of torsion.
If $\pi_1 G$ is free abelian, we may substitute the maximal torus $T$
for $G$ and since $H^*(BT)$ is torsion-free, this turns out to be automatic.
Otherwise, the hypothesis on torsion primes of $G$ and $X$
requires a more complex argument to arrive at the same conclusion.
}

\rev{
Once this form is shown, one has a description of the
associated graded ring $\Direct F_p/F_{p+1}$ of $\K(X_G;\kk)$
as a module over the associated graded ring $\Direct \ul F_p/\ul F_{p+1}$ 
of $\K(BG;\kk)$.
The first filtrand $F_1$ is the kernel of 
$\K(X_G;\kk) \lt \K(X;\kk)$,
which we would like to show is $\wt K^*(BG;\kk)\K(X_G;\kk)$.
That is, we would like to see $F_1$ is the ideal generated by 
the image of $\ul F_1 = \wt K^*(BG;\kk)$ in $F_0 = \K(X_G;\kk)$.
The description of the associated graded shows us 
that this is so ``up to an error term'' given by $F_p$,
but a lemma of Atiyah--Segal on this filtration shows that for high enough $p$
this error term also lies in the ideal generated by $\wt K^*(BG;\kk)$,
concluding the proof.
}
\begin{proof}
  \dead{If $\pi_1 G$ is free abelian,
	consider a maximal torus $T$ of $G$.
	The forgetful map, assumed surjective, 
	factors through $\KT(X;k) \lt \K(X;k)$,
	which must then be surjective as well.
	By \Cref{thm:HLGT}, the $k$-\WEFity of the restricted $T$-action
	would imply that of the original $G$-action,
	so it will be enough to show the former.
	Thus if $\pi_1 G$ is free abelian, we may assume $G$ is a torus.}
	
	We adopt the notation of \Cref{thm:KK-to-K},
	in which the augmentation ideal of $R(G;\kk)$ is $\defm I = I(G;k)$,
	and its completion $\wh I = \wt K^*(BG;k)$
	is annihilated by the map
	$\smash{\defm{\wh f}\: \K(X_G;\kk) \lt \K(X;\kk)}$ of completions.
	From \Cref{thm:KK-to-K}\ref{thm:KK-K-surj},
	we know $\wh f$ is surjective, so 
	the induced  $\defm{\bar{\wh f}}\mn\:\K(X_G;\kk)/\wh I \. \K(X_G;\kk) \lt \K(X;\kk)$ is also.
	We will show that $\ker \wh f = \wh I \. \K(X_G;\kk)$,
	so that $\bar{\wh f}$ is injective as well,
	and hence by \Cref{thm:KK-to-K}\ref{thm:KK-K-surj},
	the reduction $\bar f\:\K(X_G;\kk)/\wh I \. \K(X_G;\kk) \lt \K(X;\kk)$,
	surjective by hypothesis, is also injective,
	giving $\kk$--weak \eqfity.
	
	The $\K(BG;k)$-algebra structure on $\rev{\K(}X_G\rev{;\kk)}$ \num{321}
	under consideration  
	is that induced from the Borel fibration
	$\smash{X \to X_G \os{\defm\xi} \to BG}$.
	Although the total space $X_G$ is noncompact,
	its \AHLSSS $\defm {E_\bl(\xi)}$ 
	converges to $\K(X_G;k)$ by \Cref{thm:AHLSSSML} 
	because $\big(\K(X_{p,G};k)\mnn\big)$ satisfies the Mitt\"ag--Leffler 
	criterion~\cite[Cor.~2.4]{atiyahsegalcompletion},
	where $\defm{X_{p, G}}$ is $\xi\- B_p G$
	for $\defm{B_p G}$ the $p$-skeleton in
	a CW-structure on $BG$ \cite[Thm.~5.1]{milnor1956universalII}
	such that 
	$B_0 G = {*}$ is the basepoint. 
	%
	As $\pi_1 BG = \pi_0 G = 0$,
	we have $E_2(\xi) = \H\big(B G; \K (X;k)\mnn\big)$,
	and since the surjective fiber restriction $\smash{\wh f}$
	can be written as the composite 
	\[
	\K(X_G;k) 
	\longepi 
	\Ei\col(\xi)
	\longinc 
	E_2\col(\xi) 
	= 
	H^0\big(B G;\K (X;k)\mnn\big) 
	\isoto 
	\K(X;k)
	\]
	we see the inclusion
	$\Ei\col(\xi)
	\longinc
	E_2\col(\xi)$ is also surjective,
	meaning all  
	differentials originating in the left columns $E_\bl\col(\xi)$ vanish.

	%

	%

	If we view the horizontal maps in the square
	\[
	\xymatrix@C=1.75em{
		X_{G} 	\ar[r]^{{\xi}}\ar[d]_{\xi}& B G \ar[d]^\id\\
		B G 	\ar[r]_\id		& B G
	}
	\]
	as a bundle map, we may consider the associated map of \AHLSSS{s}.
	That of the right bundle is just the \AHSS
	$E_\bl(\id_{B G})$
	starting with $E_2(\id_{B G}) = H^*(B G;\kk)$ 
	and converging to $\K(B G;\kk)$ by the case $X = {*}$ of
	the convergence argument above.
	Because the differentials on the left column vanish,
	we can hope the bundle map $(\xi,\id_{B G})$
	induces an isomorphism of spectral sequences
	$E_\bl(\id_{BG}) \ox_k \K(X;k) \isoto E_\bl(\xi)$.
	
	\rev{If $\pi_1 G$ is free abelian,
	consider a maximal torus $T$ of $G$.
	The forgetful map, assumed surjective, 
	factors through $\KT(X;k) \lt \K(X;k)$,
	which must then be surjective as well.
	By \Cref{thm:HLGT}, the $k$-\WEFity of the restricted $T$-action
	would imply that of the original $G$-action,
	so it will be enough to show the former.
	Thus if $\pi_1 G$ is free abelian, we may assume $G$ is a torus $T$.}
	If \dead{$G = T$ is a torus and} $k = \Z$, then $E_2(\id_{B T}) = H^*(BT;\Z)$
	is torsion-free, so $E_2(\xi) \iso H^*(B T;\Z) \ox \K X$.
	As $E_\bl(\id_{B T})$ collapses~\cite[Thm.~2.4]{atiyahhirzebruch}
	and the map ${E_\bl(\id_{BT}) \lt E_\bl(\xi)}$
	of spectral sequences
	preserves differentials, 
	it follows the differentials on the image vanish.
	Since the differentials vanish on the image of $\K(X;\Z)$
	as well and are derivations,
	$E_\bl(\xi)$ collapses, so the isomorphism
	of $E_2$ pages persists to an isomorphism of $E_\infty$ pages.
	
	If $\rev{\pi_1}G$ is not \dead{a torus} \rev{free abelian},
	we have to worry about torsion
	before constructing a spectral sequence map.
	Since $\Z \longinc k$ is a localization, 
	it is a flat extension,
	so $H^*(BG;k) \iso H^*(BG;\Z) \ox k$,
	and also $\K(X;k) \iso \K(X;\Z) \ox k$ by definition.
	As the torsion primes of
	$H^*(BG;\Z)$ and $H^*(G;\Z)$ are equal 
	(consider the \SSS of the universal bundle),
	our hypothesis on $k$ ensures that 
	$H^*(BG;k)$ and $H^*(X;k)$ share no torsion primes.
	Moreover, the \AHSS $E_\bl(\id_{BG}) \ox \Q$ collapses,
	so the image of each differential 
	$\defm{\ul d_r}$ in $E_\bl(\id_{BG})$ 
	is torsion
	and it follows~\cite{MO:AHSStorsion} 
	the torsion submodules of $E_r(\id_{B G})$
	are weakly decreasing in $r$,
	so the sets of torsion primes for each module 
	are weakly decreasing as well.
	By the compactness assumption on $X$,
	the abelian groups $E_\bl^{p,0}(\id_{BG})$ and $K^q(X;\Z)$ 
	are finitely generated
	for all $p$ and~$q$,
	so $\Tor^{-1}_\Z\big(E_\bl^{p,0}(\id_{BG}),K^q(X;\Z)\mnn\big)$
	decomposes as a direct sum of cyclic groups 
	$\Tor(\Z/m,\Z/n) \iso \Z/(m,n)$,
	where $m$ and $n$ are prime powers. 
	Particularly, each the order of each nonzero summand is a power of 
	a common torsion prime $\ell$ of $H^p(BG;\Z)$ and $K^q(X;\Z)$.
	Since $k$ is a flat extension of $\Z$,
	the fact all these primes are invertible in $k$ ensures each
	$\smash{\Tor^{-1}_\Z\big(E_\bl^{p,0}(\id_{BG}),K^q(X;k)\mnn\big)}$ 
	vanishes and also that 
	$\smash{\Tor^{-1}_\Z\mn\big(\mn\im \ul d_r^{p,0},K^q(X;k)\mnn\big)}$ 
	is zero.
	Then by the ``all-cohomology'' version of the 
	universal coefficient theorem~\cite[Thm. 5.5.10, p.~246]{spanier},
	since $k$ is a principal ideal domain and $K^*(X;k)$ is of finite type
	we have $E_2(\xi) = \H\big(BG;\K(X;k)\mnn\big) 
	\iso H^*(B G;k) \ox_k \K(X;k)$.
	By assumption, $d_2$ vanishes on the $\K(X;k)$ factor of $E_2(\xi)$,
	so it can be identified with $\ul d_2 \ox \id_{\K(X;k)}$.
	Inductively, tensoring the short exact sequences
	$\ker \ul d_r \to E_r(\id_{BG}) \to \im \ul d_r$
	and $\im \ul d_r \to \ker \ul d_r \to E_{r+1}(\id_{BG})$
	with $\K(X;k)$ over $k$ and taking cohomology,
	the vanishing of the Tors we have troubled over yields
	identifications $E_r(\xi) \iso E_r(\id_{BG}) \ox_k \K(X;k)$
	and $d_r = \ul d_r \ox \id_{\K(X;k)}$
	for all $r$.
	
	Thus, whether or not $\rev{\pi_1}G$ is \rev{free abelian} \dead{a torus},
	comparing $\Ei$ pages gives an isomorphism
	\quation{\label{eq:gr-total}
		\gr \K(X_G;k) \iso \big({\gr\K(B G;\kk)}\big) \ox_\kk \K(X;\kk)
		\mathrlap,
	}
	where
	the filtrations of $\K(BG;k)$ and $\K(X_G;k)$ 
	defining the associated graded modules
	are given by
	\eqn{
		\defm{\ul F_p{}^{p+q}} &\ceq\smash{
			\im\mn\big(K^{p+q}(BG,B_{p-1} G;k) 	
			\lt 
			K^{p+q}(BG;k)\mnn\big)}\mathrlap,\\
		\defm{F_p{}^{p+q}} &\ceq 
		\smash{\im\mn\big(K^{p+q}(X_G,X_{p-1,G};k) 
			\lt
			K^{p+q}(X_G;k)\mnn\big)}
	}
	respectively, and 
	we may as well consider only $\defm{\ul F_p} \ceq \ul F_p{}^0$
	since $K^1 BG = 0$.
	Since we took $X_{p,G} = \xi\-B_p G$ in the topological filtration,
	we have $\xi^* \ul F_p \leq F_p{}^{0}$,
	so the $\K(BG;\kk)$-module structure on $\K(X_G;\kk)$
	induced by $\xi$
	satisfies
	$\ul F_{p'} \. F_{p}{}^{p+q} \leq F_{p'+p}{}^{p+q}$.
	There is thus an action of ${\gr \K(BG;\kk)}$
	on $\gr \K(X_G;\kk)$
	making \eqref{eq:gr-total}
	an isomorphism of $\big({\gr \K(BG;\kk)}\mn\big)$-modules.

	Writing the isomorphism \eqref{eq:gr-total}
	in terms of the summed filtrands 
	$\defm{F_p} \ceq \Direct_q F_p{}^{p+q}$,
	we have 
	\[
	F_p/F_{p+1} 
	\iso 
	\xt\kk{(\ul F_p/\ul F_{p+1})}{\K X}
	\mathrlap.
	\]
	Using the fact that \eqref{eq:gr-total}
	is a $\big({\gr \K(BG;\kk)}\big)$-module isomorphism%
	\rev{ and the observation that $F_0 = K^*(X_G;\kk)$},
	we then extract the relations
	\[
	F_p =  \ul F_p \. F_0 + F_{p+1} 
	\leq 
	\ul F_1 \. F_0 + F_{p+1}
	\mathrlap.
	\]
	By inductive substitution,
	one finds for arbitrarily high $p$ that
	\[
	F_1 = \ul F_1 \. F_0 + F_p\mathrlap.
	\]
	As $B_0 G$ was chosen to be the basepoint ${*}$,  
	we have $\ul F_1 = \ker\mn\big(K^0(BG;\kk) \to K^0({*};\kk)\mnn\big) = \smash{\IDEAL}$.
	The $F_\bl$-topology on $\smash{\K(X_G;k)}$ is the same as
	the $\smash{\IDEAL}$-adic topology~\cite[Cor.~2.3]{atiyahsegalcompletion},
	so some $F_p$ is contained in ${\IDEAL}\.\K(X_G;k) \rev{{}= \ul F_1\.F_0}$,
	meaning $\ker\smash{\wh f} = F_1 \rev{{} = \ul F_1\.F_0}$\dead{is generated by $\smash{\IDEAL}$}.
	But as $\xi\- B_0 G = X$
	is the fiber of the Borel fibration over the basepoint, 
	we may identify the fiber restriction 
	$\smash{\wh f}$ with $F_0 \longto F_0/F_1$,
	whose kernel $F_1$ we have just seen is \rev{$\ul F_1\.F_0$}\dead{generated by $\smash{\IDEAL}$}.
	\rev{Thus $\ker \wh f = \wh I \. \K(X_G;\kk)$, as was to be shown.}
	\epf
	
	\begin{remarks}
		\emph{(a)}
		Harada and Landweber chose the name \emph{weak equivariant formality}
		because equivariant formality in \EC is equivalent to the collapse of the
		\SSS of the Borel fibration $X \to X_G \to EG$,
		which implies the ``strong''
		condition that $H_G^* X \iso H_G^* \ox_\Q H^*X$ as an $H_G^*$-module.
		This is equivalent to the collapse of the \AHLSSS 
		with rational coefficients, by the Chern character isomorphism,
		but by \Cref{ex:unfree}, we cannot say anything similarly strong for $\KG(X;k)$,
		at least without the assumption on $\pi_1 G$.
		We can however get a similar result for 
		$\KG(X;k)\compl = \K(X_G;k)$ by a Leray--Hirsch--style argument 
		for coefficient rings $k \leq \Q$ such that $\K(X;k)$ is a free a $k$-module.
		
		\smallskip
		
		\emph{(b)} Harada and Landweber note that \WEFity 
		is equivalent to the edge map $\Z \ox_{RG} \KG X \lt \K X$ 
		of the \HKSS $\Tor_{RG}(\Z,\KG X) \implies \K X$ being an isomorphism,
		and note that if $\pi_1 G$ is torsion-free,
		one can certainly guarantee this (and collapse) by assuming $\Tor^{\leq -1}$ vanishes%
		~\cite[Prop.~4.4]{haradalandweber2007}.
		The combination of \Cref{thm:weak} and any instance of \Cref{thm:main}
		with $\rk H < \rk G$ (e.g., \Cref{thm:TNHZ,thm:gensymm})
		shows this sufficient condition is not necessary.
		
		\smallskip
		
	\end{remarks}
	
	\bex\label{ex:unfree}
	One is tempted to analogize from Borel cohomology that 
	$\KG(X;\kk) \lt \K(X;\kk)$ surjective implies $\KG(X;\kk)$ is free over $R(G;\kk)$.
	This would certainly simplify the proof of \Cref{thm:weak}.
	One reason not to yield to this temptation is as follows.
	Let $T$ be a maximal torus of $G = \PSU(3)$ and consider the transitive action 
	of $G$ on $G/T$.
	The induced map $\HT \lt H^*(G/T)$
	is surjective with image $\HT \ox_{H_G^*} \Q$
	by \Cref{thm:formalCDGA} because $\rk T = \rk G$,
	but by \Cref{thm:equalrank},
	this is $RT \ox_{RG} \Q \iso \K(G/T;\Q)$,
	so it follows $\KG(G/T;\Q) \iso R(T;\Q) \to \K(G/T;\Q)$
	is surjective. 
	But $RT \ox \Q$ is not free over $R\PSU(3) \ox \Q$\rev{; 
	see Brylinski--Zhang~\cite[\S7]{brylinskizhang}\num{Ex.2.27}
	for a description of $R\PSU(3)$ as a subring
	of $RT$, and Fok~\cite{MO:Fok} 
	for an explicit statement of the algebra generators $X,Y,Z$ 
	(the $a,b,c$ of \Cref{rmk:completion-weak})}. 
	\dead(Note \dead{however} that $\a\: RT \lt \K\mnn\big(\PSU(3)/T\big)$ is 
	\emph{not} surjective integrally%
	~\cite[Thm.~6]{matsunagaminami}.\dead)
	\eex

\section{The isotropy action}\label{sec:main}
In this section we obtain \Cref{thm:main} and the weaker but more general 
\Cref{thm:condcollapse}.

\begin{definition}	
	If $\defm H$ is a closed, connected subgroup of the compact, connected Lie group $G$---as it always will be---then we call
	$(G,H)$ a \defd{\ccpair}.
	The left multiplication action of $H$ on $G/H$ is the \defd{isotropy action}.
\end{definition}
\begin{definition}
Given a graded $\kk$-algebra $A$ and graded 
$A$-modules $M$ and $N$,
we will write $\defm{\Tor^{-p}_A(M,N)} = \Tor^{A}_p(M,N)$
and $\defm{\Tor_A(M,N)} = \Direct_{p \geq 0} \Tor^{-p}_A(M,N)$.
Such a module is bigraded, with the other degree $\defm q = i+j$ 
coming from the degree $\dead i \rev j$ in $N$ and the internal degree $i$
of a projective $A$-module resolution $P^{p,i}$ of $M$.
The motivation for this grading
is that if 
$M$ and $N$ are graded $A$-modules,
one of which is flat, 
equipped with $A$-module differentials of degree $+1$,
then $n = -p+q$ is the total degree in an 
algebraic K\"unneth spectral sequence 
$\Tor_{A}^{-p,q}(M,N) \implies H^{-p+q}(M \ox_A N)$.
The \EMSS and the other \KSS we will discuss 
share this bigrading.\num{p13}
We notationally suppress \dead{these} 
both these gradings whenever practicable.
\end{definition}

\begin{discussionsilent}
	The principal tools for computing the related \BEC ring
	$\HH(G/H)$ are the \SSS 
	of a fibration $G \to G_{H \x H} \to BH \x BH$
	and the closely related \EMSS of the homotopy pullback square on 
	the left in the following figure---
	
	\quation{\label{eq:Gfibersquare}
		\begin{gathered}
			\xymatrix@R=2.5em@C=1em{
				G_{H \x H}\ar[r]\ar[d]		& BH\ar[d] 		&&&&&
				(X \x Y)_G 	\ar[r]\ar[d]	& Y_G \ar[d]\\
				BH \ar[r]  					& BG 			&&&&& 	
				X_G			\ar[r]			& BG
			}
		\end{gathered}
	}
	
	\medskip
	
	\noindent---which is the special case $X = Y = G/H$
	of the square on the right.
	The \EMSS of the right square is the general \KSS in \BEC,
	and in our case of interest the sequence collapses; 
	in fact, $\H_H(G/H;\Q)$ is isomorphic \emph{as a ring} to the $E_2$ page.
	
	The K-theoretic analogue is a \KSS 
	\rev{due to Hodgkin~\cite{hodgkin1975kunneth},}
	with intended target $\KG(X \x Y)$\dead, 
	\dead{due to Hodgkin~\cite{hodgkin1975kunneth},}
	\rev{and}
	beginning at $E_2 = \Tss_{RG}(\KG X, \KG Y)$.
	This sequence strongly converges, 
	but not usually to the intended target~\cite[Thm.~5.1]{hodgkin1975kunneth}.
	It does converge to the desired $\KG(X \x Y)$,
	however, 
	if $G$ is compact and connected with free abelian fundamental group%
	~\cite{hodgkin1975kunneth,snaith1972kunneth,mcleod1979kunneth}
	\rev{(see \Cref{rmk:convergence} for an account of the proof)}.
	When $X = Y = G/H$, one applies the standard equivalences
	\[\KH(G/H) = \K_{H}(\xt G G {G/H}) \iso 
	\K_{G \x H}(G \x G/H) \iso 
	\K_G\big((G/H) \x (G/H)\mnn\big)\]
	to see $\KH(G/H)$ is the target of the \HKSS
	beginning with the $E_2$ page 
	$\Tor_{RG}\big(\KG(G/H),(\KG(G/H)\mnn\big) \iso 
	\Tor_{RG}(RH,RH)$.
	To obtain a\dead{n}\num{p13} collapse result analogous 
	to the one for the \EMSS, 
	we compare the two spectral sequences\rev:
\end{discussionsilent}

\SPECMAP

As so often with spectral sequence constructions,
the proof of Theorem \rev{3.5} 
\dead{4.4}\num{p13}
is logically almost independent from the use we make of it,
so we black-box it here and take up the proof in \Cref{sec:mapkunneth}.
Our intended use is to give a sufficient condition for the 
collapse of the \HKSS of interest.

When $X = Y = G/H$,
the map of $E_2$ pages is
\[
\Tor_{RG}(RH,RH) \lt \Tor_{H_G^{**}}(\HHp,\HHp)
\]
and that on $\Ei$ pages is
\[
\gr \KH(G/H) \lt \gr \HHp(G/H).
\]
As noted in \eqref{eq:extension},
these maps are induced by rational completion 
$M \mapsto M \ox \Q \mapsto (M \ox \Q)\compl$
with respect to the augmentation ideal $IG \ox \Q$.


\bprop\label{thm:condcollapse}
Let $k$ be a subring of $\Q$.
If $\defm T = \Tor_{R(G;\kk)}\!\big(R(H;\kk),R(\dead{K}\rev{H};\kk)\mnn\big)$
is such that $T \lt T \ox_{R(G;\kk)} H_G^{**}$
is an injection,
then the \KSS $T \implies \KH(G/H;k)$ collapses.
\eprop
\bpf
Recall that $\HHp$ is the completion of $R(H;\Q)$ with respect to 
$I(H;\Q)$ but also with respect to $I(G;\Q)$, and similarly for $H^{**}_{\dead{K}\rev{H}}$~\cite[Prop.~3.9]{segal1968representation}.
Thus if $\defm{P^\bl}$ is a projective $R(G;\kk)$-module resolution of $R(H;\kk)$,
then $P^\bl \ox_{R(G;\kk)} H_G^{**}$
is a projective $H_G^{**}$-module resolution of $R(H;\kk) \ox_{R(G;\kk)} H_G^{**} = \HHp$.
Since $R(G;\kk)$ and $R(\dead{K}\rev{H};\kk)$ are Noetherian~\cite[Cor.~3.3]{segal1968representation}
and $T$ is finitely generated over $R(\dead{K}\rev{H};\kk)$ 
and hence~\cite[Prop.~3.2]{segal1968representation} over $R(G;\kk)$, 
the completion $(T \ox \Q)\compl$ with respect to $I(\dead{K}\rev{H};\Q)$ (or equivalently, $I(G;\Q)$)
is isomorphic to $T \ox_{R(G;\kk)} H_G^{**}$ and $ T\ox_{R(\dead{K}\rev{H};\kk)} H_G^{**}$~\cite[Prop.~10.13]{atiyahmacdonald}.
Noting also that $H^{**}_{\dead{K}\rev{H}}$ is flat over $R({\dead{K}\rev{H}};\kk)$ by \Cref{thm:extflat},
so that $- \ox_{R({\dead{K}\rev{H}};\kk)} H^{**}_{\dead{K}\rev{H}}$ commutes with cohomology,
we see $\Tor_{H_G^{**}}(\HHp,H^{**}_{\dead{K}\rev{H}})$ can be computed as
\eqn{
H^*(P^\bl \!\!\ox_{R(G;\kk)}\! H_G^{**} \!\ox_{H_G^{**}}\! H^{**}_{\dead{K}\rev{H}})
	\,\iso\,
H^*(P^\bl \!\!\ox_{R(G;\kk)}\! H^{**}_{\dead{K}\rev{H}})	
&	\,\iso\,
\H\big(P^\bl\!\! \ox_{R(G;\kk)}\! R({\dead{K}\rev{H}};\k) \!\! \ox_{R({\dead{K}\rev{H}};\kk)}\! H^{**}_{\dead{K}\rev{H}}\big)
\\&	\,\iso\,
\H\big(P^\bl \ox_{R(G;\kk)} R({\dead{K}\rev{H}};\k)\mnn\big) \ox_{R({\dead{K}\rev{H}};\kk)} H^{**}_{\dead{K}\rev{H}}
	\,=\,
T \ox_{R({\dead{K}\rev{H}};\kk)} H^{**}_{\dead{K}\rev{H}}
\mathrlap,
}
%
so the map we are assuming injective 
is the map of $E_2$ pages in the 
map of spectral sequences we construct in \Cref{thm:eqkunnethchern}.
As $H_G^{**}$ is flat over $R(G;\kk)$ and
the differentials are $R(G;\kk)$-linear,
extension by $H_G^{**}$ commutes with taking cohomology,
so by induction, the map on each subsequent page is also
$E_r \lt E_r \ox_{R(G;\kk)} H_G^{**}$.
By assumption, the $E_2$ map is an injective cochain map,
so as the target sequence collapses, 
the differential $d_2$ is trivial
and $E_2 = E_3$. 
But one can argue the same of $d_3$,
and by induction the source sequence collapses as well.
\epf

\bcor
Let $(G,H)$ be a pair of compact, connected Lie groups such that
$\pi_1 G$ is free abelian.
Then $H$ acts equivariantly formally on $G/H$ if and only if
the map 	
\[
\Tor_{\id}(\e,\id)\:
\Tor_{R(G;\Q)}\!\big(R(H;\Q),R(H;\Q)\mnn\big)
\lt 
\Tor_{R(G;\Q)}\!\big(\Q,R(H;\Q)\mnn\big)
\]
is surjective.
\ecor
\begin{proof}
	Since the $IG$-adic and $I\dead{K}\rev{H}$-adic topologies on $R\dead{K}\rev{H}$ agree~%
	\cite[Prop.~3.9]{segal1968representation},
	applying the flat extension $- \ox_{RG} H_G^{**}$
	to $\Tor_{\id}(\e,\id)$ 
	gives \[\Tor_{H_G^{**}}(\HHp,\HHp) \lt \Tor_{H_G^{**}}(\Q,\HHp)\mathrlap.\]
	By the collapse of the \EMSS of the square 
	on the left in \eqref{eq:Gfibersquare},
	this is a way of writing $\HHp(G/H) \lt H^*(G/H)$,
	yielding the forward implication.
	For the backward, note that under the hypothesis on $G$,
	the map in question
	is the $E_2$ page of a map of {\HKSS}s 
	converging to the map $\KH(G/H;\Q) \lt \K(G/H;\Q)$.
	The sequence converging to $\K(G/H;\Q) \iso \Tor_{R(G;\Q)}\!\big(\Q,R(H;\Q)\mnn\big)$ always collapses,
	so since the $\Ei$ page of the \dead{other} \rev{domain} sequence is a 
	subquotient of the $E_2$ page, 
	if $\KH(G/H;\Q) \lt \K(G/H;\Q)$ is surjective,
	so \emph{a~fortiori} is $\Tor_{\id}(\e,\id)$. 
\end{proof}

\brmk
\Cref{thm:condcollapse} should not come as such a surprise.
Minami proved \cite[Thm.~2.1]{minami1975symmetric}
that if $(G,\dead{K}\rev{H})$ is a \ccpair with $\pi_1 G$ torsion-free,
then the K{\"u}nneth spectral sequence converging to $\K(G/\dead{K}\rev{H})$
collapses at $E_2 = \Tor^{*,*}_{RG}(\Z,R\dead{K}\rev{H})$
and pursued this~\cite{minami1975symmetric,minami1976symmetric} 
to a complete description of the K-theory of symmetric spaces
$G/\dead{K}\rev{H}$ with $\pi_1 G = 0$.
\ermk

If the Tor is well-behaved, we can obtain the collapse
from first principles without direct reference 
to cohomology and resolve the multiplicative
extension problem manually. 
This path will lead us to \Cref{thm:main}.

\begin{lemma}\label{thm:collapse-trivial}
	Let $G$ be compact, connected Lie group 
	with torsion-free fundamental group,
	\dead{and} $X$ and $Y$ 
	compact $G$-spaces\rev, and $\kk$ a coefficient 
	ring such that $\KG(X;\kk)$ and $\KG(Y;\kk)$
	are finite $R(G;\kk)$-modules.
	Suppose further that $\Tor_{R(G;\kk)}\!\big(\KG(X;\kk),\KG(Y;\kk)\mnn\big)$
	is generated as an algebra over $\Tor^0$ by generators in $\Tor^{-1}$
	and that $\Tor^0$ contains no $2$-torsion.
	Then the \HKSS converging to $\KG(X \x Y;\kk)$ collapses at $E_2$.
\end{lemma}
\begin{proof}
$\Tor = \Direct \Tor^{p,q}$ is the $E_2$ page of the \HKSS, and is generated as an algebra
by elements in columns $p = 0$ and $p = -1$.
Since the differential $d_2$ increases $p$ by $2$,
it vanishes on all generators, and since it is a derivation,
it vanishes entirely. By induction, we have $\Tor = E_2 = E_\infty$.
\end{proof}
\begin{lemma}\label{thm:collapse-trivial-alg}
	Assume the conditions of \Cref{thm:collapse-trivial}
	and moreover that 
	$\Tor$
	is the \emph{free commutative graded} algebra over $\Tor^0$ 
	on a sequence $\vec z$ of \dead{homogeneous} generators in $\Tor^{-1}$
	\dead{and that $\Tor^0$ contains no $2$-torsion}.
	Then 
	\[
	\KG(X \x Y;\kk) \iso \Tor_{R(G;\kk}\!\big(\KG(X;\kk),\KG(Y;\kk)\mnn\big)
	\]
	as graded $\kk$-algebras.
\end{lemma}
\begin{proof}
We already have collapse,
so it remains to resolve extension problems.
The filtration $(\defm{F_p})_{p \in \Z}$
of $\KG(X \x Y;\kk)$
for which $E_\infty = \Tor$ is the associated graded
is nonpositive and particularly has filtrand $F_1 = 0$,
so that $\Tor^0 = F_0/F_1 = F_0$
is a naturally identified as a subring of $\KG(X \x Y;\kk)$.
The map witnessing this is the map $\l$
discussed in \Cref{ex:SO(3)} and \Cref{rmk:convergence},
induced by the exterior product
  $\KG\rev(X\rev{;\kk)} \dead{\x}\rev{\ox_{R(G;\kk)}} \KG \rev(Y\rev{;\kk)} \lt \KG(X \x Y\rev{;\kk})$.\num{477}

We now seek to identify a complementary tensor factor 
$\smash{\ext \wP}$ of $\KG(X \x Y;\kk)$
lifting the factor $\ext[\vec z]$ of $\Tor$.
For this, lift each exterior generator $z_i$
to an element $\defm{\wt z_i}$ of $F_{-1}K^1_G(X \x Y;\kk)$
and each polynomial generator $z_j$ to an element
$\defm{\wt z_j}$ of $F_{-1}K^0_G(X \x Y;\kk)$
and write $\defm \wP$ for the $\kk$-span of the $\wt z_i$ and $\wt z_j$.

Then, to see that the $(\im \l)$-algebra
generated by $\wP$
is all of $\KG(X \x Y;\kk)$,
note that the associated tensor factor of 
$\gr \KG(X \x Y;\kk)$ generated by $\vec z$
spans it as a module over $\Tor^0 \iso F_0 = \im \l$.
Thus $\smash{\wP \. \im \l}$ represents $F_{-1}/F_0$
in the associated graded 
and $\smash{\im \l + \wP \.\im \l}$ generates $F_{-1}$.
Similarly $(\kk + \wP + \cdots + \wP^p)\.\im \l$ generates $F_{-p}$ for each $p$,
and the claim follows by induction. 

We claim the lifts $\wt z_i$ square to zero.
Indeed, by graded commutativity,
$\wt z_i^2 = - \wt z_i^2 $ is $2$-torsion, lying in $F_{-2}$.
Since all $2$-torsion in $F_{-2}/F_{-1}$ and $F_{-1}/F_0$ and $F_0$
is trivial, it follows $\wt z_i^2 = 0$.
Now we see the previous induction represents $F_p$ as the direct sum
$(\kk \+ \wP \+ \cdots \+ \wP^p)\.\im \l$,
which shows $\KG(X \x Y;\kk)$ is a free $\Tor^0$-algebra,
  and in particular shows $\ext\wP \ceq \kk[\wt z_j] \ox_{\rev\kk} \ext[\wt z_i]$
  is a free \CGA with $\KG(X \x Y;\kk) \iso (\im \l) \ox_{\rev\kk}\ext\wP$.
\end{proof}

\begin{corollary}\label{thm:nicecorollary}
	Let $G$ be a compact, connected Lie group with
	torsion-free fundamental group,
	$H$ and $K$ closed, connected subgroups, 
	and $\kk$ a subring of $\Q$ such that $\Tor_{R(G;\kk)}\big(R(H;\kk),R(K;\kk)\mnn\big)$
	is an exterior algebra $\Tor^0 \ox_{\rev\kk} \ewP$ over
	$\Tor^0 = R(H;\kk) \otimes_{R(G;\kk)} R(K;\kk)$
	on elements of $\Tor^{-1}$ and the latter is $2$-torsion--free.
	Then $\KH(G/K;\kk)$ is also isomorphic to 
	$R(H;\kk) \otimes_{R(G;\kk)} R(K;\kk) \ox \ewP$
	as a $\Z/2$-graded ring.
\end{corollary}

\rev{It is not clear in general what pairs of 
	subgroups $H$ and $K$ satisfy such
	generous hypotheses,
	but the case $H = K$ of the isotropy action fortunately does
	in several of the best-studied equivariantly formal cases.}
The predecessor work~\cite{carlsonfok2018}, 
identified several equivalent characterizations of
equivariant formality of the isotropy action,
including the condition that the kernel of the restriction map 
$R(G;\Q) \lt R(K;\Q)$
descend to a vector subspace of dimension $\rk G - \rk K$
in the space $I(G;\Q)/I(G;\Q)^2$ of indecomposables.
A closely related but generically stronger condition
is that the kernel of the restriction
be generated as an ideal by a regular sequence of $\rk G - \rk K$
elements.

\MAIN*

\begin{proof}
	We show our hypotheses imply those of \Cref{thm:nicecorollary}
	with $K = H$.
	
		We have assumed $R(H;\kk)$ is flat 
		over the image $\defm R$ of $R(G;\kk) \to R(H;\kk)$,
		which is isomorphic to the quotient of $R(G;\kk)$ by the ideal
		$\defm\fa = (\rho_j)$ generated by some regular sequence 
		of 
		elements 
		$\defm{\rho_j}$.
	Let $\defm\wP = \kk\{\defm{z_j}\}$ be the free $\kk$-module on 
	as many generators,
	so that the associated Koszul complex 
	$\ewP \ox_\kk R(G;\kk)$ with differential given by
	$z_j \mapsto \rho_j \mapsto 0$ 
	is an $R(G;\kk)$-module resolution of $R$.
	By flatness we have
	\[
	\Tor_{R(G;\kk)}\big(R(H;\kk),R(H;\kk)\mnn\big)
		\iso
	R(H;\kk) \ox_R\ \Tor_{R(G;\kk)}\big(R,R(H;\kk)\mnn\big)\mathrlap.		
	\]
	For the second tensor factor, one uses the Koszul complex to compute
	\eqn{
		\Tor_{R(G;\kk)}\big(R,R(H;\kk)\mnn\big)
		&\iso
		\H\big(
		\ewP \ox_\kk R(G;\kk) \ox_{R(G;\kk)} R(H;\kk)
		 \big)
		\\&		\iso 
		\H\big(\ewP \ox_\kk R(H;\kk)\mnn\big)
			\\&		\iso 
		\ewP \ox_\kk R(H;\kk)
		\mathrlap,
	}	
	since the differential on $\ewP \ox_\kk R(H;\kk)$ is 
	determined by $z_j \mapsto \rho_j \pmod\fa = 0 \in R(G;\kk)/\fa \iso R$
	and hence is zero.
	The $R$-module structure on this Tor arises via the inclusion of $R$ in the $R(H;\kk)$
	tensor factor, so
		\[
		\Tor_{R(G;\kk)}\!\big(R(H;\kk),R(H;\kk)\mnn\big)
		\iso
		R(H;\kk) \ox_R \big(\ewP \ox_\kk R(H;\kk)\mnn\big)
		\iso
		\big(R(H;\kk) \ox_R R(H;\kk)\mnn\big) \ox_\kk \ewP
		\mathrlap,
		\]
	but $R(H;\kk) \ox_R R(H;\kk) \iso R(H;\kk) \ox_{R(G;\kk)} R(H;\kk)$
	by the definition of $R$.
	
	By definition we have $R(H;\kk) \ox_R R(H;\kk) \iso R(H;\kk) = \Tor^0$
	while $\wP$ lies in $\Tor^{-1}$. 
	We will be able to apply \Cref{thm:nicecorollary}
	once we see $\Tor^0$ 
	contains no $2$-torsion.
	But $R \ox_{R(G;\kk)} R \iso R \leq R(H;\kk)$ is a 
	\dead a free $\kk$-module, so multiplication by $c \in \kk\less\{0\}$ is injective,
	and as $R(H;\kk)$ is flat over $R$, multiplication by $c$
	is injective on $R(H;\kk) \ox_{R(G;\kk)} R(H;\kk)$ as well.

	To see $\K(G/H;\kk)$ is isomorphic to $\kk \ox_{R(G;\kk)} R(H;\kk) \ox_\kk \ewP$,
	we may apply essentially same argument with the first $R(H;\kk)$ replaced by 
	$\kk$, noting that multiplication by $2$ is injective on $\kk$ as well.
	The map of $E_2 = E_\infty$ pages is then clearly reduction modulo $I(H;\kk)$,
	and the naturality of the map~$\l$ identifies the induced map on $\Tor^0$,
	namely
	$R(H;\kk) \ox_{R(G;\kk)} R(H;\kk) \lt \Z \ox_{R(G;\kk)} R(H;\kk)$,
	with the map on the corresponding tensor factors of 
	$f\: \KH(G/H;\kk) \lt \K(G/H)$. 
	Particularly,~$f$ is surjective.
	It is not automatic that the restriction of $f$ to the exterior tensor factor $\ewP$
	can be identified with the identity\rev, because the lifts $\wt z \in F_{-1}$ 
	from $z \in \Tor^{-1} = F_{-1}/F_0$	in the proof of \Cref{thm:collapse-trivial-alg}
	were chosen separately and potentially inconsistently.
	To repair this, we may select exterior generators of $\K(G/H;\kk)$
	to be the images of generators of $\ewP \leq \KH(G/H;\kk)$ under $f$.
	
	That the number of exterior generators is $\rk G - \rk H$
	follows on comparison with cohomology,
	where we know $\HHp(G/H;\Q) \iso \HHp \ox_{H_G^{**}} \Q \ox_\kk \ewP$
	via the equivariant Chern character
	and we already know $\dim_\Q (\Q \ox_\kk \wP) = \rk G - \rk H$
	by \Cref{thm:formalCDGA}.
\end{proof}

\begin{remark}
	There is a good concrete choice of exterior generators:
	to each representation $\rho\: G\lt \Aut_\C V$ on some complex vector space $V$
	which is trivial on $H$ 
	(so that $\rho - \e(\rho)$ lies in the kernel of $RG \lt RH$), 
	Hodgkin~\cite[pp.~77--9]{hodgkin1975kunneth} 
	associates a bundle on the \rev{unreduced} suspension~$S(G/H)$, 
	given by the clutching construction
	gluing two copies of the trivial bundle~$C(G/H) \x V$
	on the cone 
	along the bounding copies of $G/H \x V$
	using the map
	$(gH,v) \lmt \big(gH,\rho(g)v\big)$.
	Since th\rev{at}\dead{is} map is equivariant with respect to the diagonal left $H$-action,
	this construction gives an element $\defm{\b_H^{G/H}}(\rho)$ of $K^1_H(G/H)$,
	and $f \o \b_H^{G/H}$ is the map Hod\rev{g}kin describes.
	
	This construction is natural in pairs $(G,H)$.
	To show the image obtained by 
	applying $\b_H^{G/H}$ to $\fb = \ker(RG \to RH)$ generates $\ewP$,
	we can look at the restriction of 
	$\defm{\b} = \b^{G/1}_1\: IG = \ker(RG \to R1) \to K^1 G$ to $\fb$,
	which is the composition
	$\fb \to K^1_H(G/H) \to K^1(G/H) \to K^1 G$
	induced by the inclusion $(G,1) \longinc (G,H)$;
	but Hodgkin~\cite{hodgkin1967lie} also showed $\b$ factors  through an injection
	from the free abelian group $IG/(IG)^2$,
	so that $\b(\fb)$ is of rank $\rk G - \rk H$.
\end{remark}

A number of cases of traditional interest readily follow.

%
%
\dead{
By
\num{549}\Cref{thm:isotfformal},
equivariant formality implies
	$G/H$ is formal. 
	Formality implies,
	by \eqref{eq:extension}, 
	that
		$\HH \ox_{H_G^*} \Q$
		is a 
		Poincar\'{e} duality algebra,
	which by standard commutative algebra~\cite[Thm.~6.2]{carlsonfok2018}
	is naturally isomorphic to 
	\quation{\label{eq:trivialcompletion}
		\xt{H_G^{**}}{\HHp}\Q 
			\ \iso\ \mnn 
		\xt{R(G;\Q)}{R(H;\Q)\mn}{\mn\xt{H_G^{**}\mnn}{H_G^{**}}{\Q}}
			\ \iso\ 
		\xt{R(G;\Q)}{R(H;\Q)\mn}{\mn\Q}\mathrlap,
	}
{
Now $R(H;\Q)$ is a polynomial ring by the assumption that
$\pi_1 H = 0$
and is the Weyl-invariant subring of $R(S;\Q)$
for $S$ a maximal torus of $H$.
Kane's variant \Cref{thm:chevalleyshepherdtodd}
of the Chevalley--Shephard--Todd theorem 
implies that $R(H;\Q)$ is free over 
$R = \im\big(R(G;\Q) \to R(H;\Q)\mnn\big)$,
giving the flatness hypothesis.
Now apply \Cref{thm:main} with $k = \Q$.
}
}


\TNHZ
\bpf
We apply \Cref{thm:main}, noting that in this case 
the Tor factor is simply $\Tor^0_{RG}(RH,RH) = \xt{RG}{RH}{RH} \iso RH$.
To recover $\K(G/H)$,
we apply \Cref{thm:weak}.
\epf

\EQUALRANK
\bpf
Let $T$ be a maximal torus of $H$ and write $\defm W$ for a Weyl group.
Restriction induces natural isomorphisms $RG \iso (RT)^{WG}$
and $RH \iso (RT)^{WH}$,
so since $WH$ is a reflection subgroup of $WG$,
it follows~\cite[Thm.~2.2]{steinberg1975pittie}
that $RH$ is free over $RG$.
Now apply \Cref{thm:main}.
To recover $\K(G/H)$,
apply \Cref{thm:weak}.
\epf

\brmk\label{rmk:maxrank}
	Although the author does not know a citation,
	this result should already be known.
	Recall McLeod's result that 
	if $\pi_1 G$ is free abelian, 
	then $\KT(G/T) \iso \K_{T\x T} G \iso {RT} \ox_{RG} {RT}$
	and the corollaries that (1) $\KT X \iso RT \ox_{RG} \KG X$ 
	and (2) if $\KT X$ is free over $RT$, then even if $\pi_1 G$ has torsion,
	$\KG X \iso \KT(X)^W$%
	~\cite[{\SS}3]{snaith1972kunneth}\cite[Thm.~4.4]{mcleod1979kunneth}.
	One finds first

\medskip

	\[
		\KH(G/T)
			\,\iso\,
		\KT(G/T)^{WH}
			\mnn\,\iso\,
		(\xt{RG}{RT}{RT})^{WH}
			\mnn\,\iso\,
		\xt{RG}{RH}{RT}
			\,\iso\,
		\KT(G/H),
	\]
and since this is also free over $RT$,
	
	\[
		\KH(G/H) 
		\mnn	\,\iso\, \mnn
		\KT(G/H)^{WH}
		\mnn	\,\iso\, \mnn
		(\xt{RG}{RT}{RH})^{WH}
		\mnn	\,\iso\,\mnn
		\xt{RG}{RH}{RH}.
	\]
\ermk

\GENSYMM
\bpf
\rev{
We reduce using some classical but appealingly elementary structure theory~\cite[\S5]{wolf1968homogeneous1}.
First we handle the case where the fundamental group of the identity component $H^0$
of $H$ is free abelian and $G$ is not necessarily semisimple.
The kernel of the canonical projection from $G$ to its abelianization is 
its commutator subgroup $G'$. 
The elements of $G/G'$ are represented in $G$ by the central torus $Z = Z(G)^0$,
which meets $G'$ in 
a finite subgroup $F$.
All these subgroups are preseved by any automorphism $\s$.
It is known that $(G')^{\ang{\s}}$ is connected~\cite[Rmk.~1.2]{stepien2002formal}, 
so $H$ is the product of it and some subgroup $Y$ which is a 
union of components of 
$Z^{\ang\s}$.
Hodgkin~\cite[Prop.~11.1]{hodgkin1975kunneth} 
showed the short exact sequence $1 \to G' \to G \to G/G' \iso Z/F \to 0$
induces a decomposition $RG \iso RG' \ox R(Z/F)$
when $G$ is connected with $\pi_1(G)$ free abelian,
but the proof works more generally if $G$ can be expressed as $A\.G'$
where $A$ is compact abelian and meets $G'$ in a central subgroup.
The map $R(Z/F) \lt RG$ is induced naturally by the projection,
whereas the map $RG' \lt RG$ is unnaturally obtained
by noting that $RG'$ is a polynomial ring on generators $\rho \in IG$
and for each $\rho$
finding $\upsilon \in IZ$ such that $\upsilon \ox \rho \in RZ \ox RG' = R(Z \ox G')$
restricts to the trivial representation
on the subgroup $\big\{(f,f\-):f \in F\big\}$,
and hence descends to a representation of $G = Z \ox_F G'$.
When $\pi_1(G^{\ang\s})$
is torsion-free, we can describe the restriction $RG \lt RG^{\ang\s}$
in terms of this decomposition and the analogous
decomposition $RG^{\ang\s} \iso R(\sfrac{Y}{\,Y\inter F}) \ox R(G')^{\ang\s}$
despite their noncanonicity.
Since $\sfrac{Y}{\,Y\inter F}$ is a subgroup of the torus $Z/F$,
the induced map $R(Z/F) \lt R(\sfrac{Y}{\,Y \inter F})$
is a surjection,
and hence, since the lifings of generators $\tau$ of $R(G')^{\ang\s}$
are of the form $\upsilon \ox \tau \in R(\sfrac{Y}{\,Y\inter F}) \ox R(G')^{\ang\s}$,
the image of $RG \lt RG^{\ang\s}$
can be identified as the tensor product of $R(\sfrac{Y}{\,Y\inter F})$
and the image of $RG'\lt R(G')^{\ang\s}$
independently of the choice of splitting.
This reduces us to study of the latter map.
}

\rev{
A simply-connected compact Lie group $G'$ 
is a direct product of simply-connected simple groups.
The image of any simple factor under $\s$ is another such factor,
so the set of simple subgroups decomposes into finitely 
many $\ang{\s}$-orbits, say those represented by the simple groups~$G_\ell$.
One can thus express $G'$ as the product
$\prod O_\ell$ of the $\ang\s$-invariant subgroups
$O_\ell = \prod_{i=0}^{|\s|-1}\s^i G_\ell$.
Evidently $(G')^{\ang\s}$ decomposes as $\prod O_\ell^{\ang\s}$.
Up to conjugacy by automorphisms, $O_\ell$
can be seen as a direct power $\prod_{i=1}^m G_\ell^{(i)}$ where
$G_\ell^{(i)} = G_\ell$ for $i \in [1,m]$
for some $m$ dividing $|\s|$, 
with $\s$ cyclically permuting the factors $G_\ell^{(i)}$ in
such a way that the restriction $G_\ell^{(i)} \lt G_\ell^{(i+1)}$ for $i < m$
is $\id_{G_\ell}$ and the restriction $G_\ell^{(m)} \lt G_\ell^{(1)}$
is a potentially nontrivial automorphism $\phi_\ell$ of~$G_\ell$.
Thus one has for $O_\ell^{\ang\s}$ the diagonal subgroup $\D G_\ell^{\ang{\phi_\ell}}$.
The inclusion $O_\ell^{\ang\s} \longinc O_\ell$
then factors as $G_\ell^{\ang{\phi_\ell}} \simto \D G_\ell^{\ang{\phi_\ell}} \inc \D G_\ell \inc G_\ell^m$,
so that $RO_\ell \lt RO_\ell^{\ang\s}$ factors as
$R(G_\ell^m) \epi RG_\ell \to R(G_\ell^{\ang{\phi_\ell}})$
Since $G_\ell$ is simply-connected and hence $RG_\ell$
is polynomial, the kernel of the first factor is a regular sequence
$\rho - \s^j \rho$ for $j \in [2,m]$, where $\rho$ runs over a set of polynomial generators.
If we can show for each $\ell$ that the kernel of the second factor is also generated by a regular
sequence and that $R(G^{\ang{\phi_\ell}})$ is free over its image, then,
we will have verified the hypotheses of \Cref{thm:main}.
}

\rev{
Thus we will consider the case $G$ is simple.
Then we know the quotient group $(\Aut G) / (\Inn G)$ of the automorphism group of $G$
by the inner automorphisms is isomorphic to the group of graph automorphisms
of the Dynkin diagram of $G$,
so that $\s$ can be expressed as the composite 
of a diagram automorphism $\tau$ and conjugation by an element of $G$.
Replacing $\s$ by some conjugate $c_g \o \s \o c_{g}\-$ and hence
$G^{\ang{\s}}$ by $c_g G^{\ang\s}$,
we may in fact~\cite[\S4]{goertschesnoshari2016} write $\s = \tau \o c_s$
for $s$ lying in any given maximal torus $S$ of $K = G^{\ang\tau}$;
this $S$ is then also a maximal torus of $H = G^{\ang\s}$.\footnote{\ 
	This sort of reduction to a distinguished pair $(G,K)$ coming
	from a diagram automorphism plays a role both in 
	St{\k{e}}pie{\'n}'s proof that generalized symmetric spaces are formal in
	the sense of rational homotopy theory~\cite{stepien2002formal}
	and Goertsches--Noshari's proof that the isotropy action of $K$ on $G/K$
	is equivariantly formal~\cite{goertschesnoshari2016}.
}
The representation rings of $K$ and $H$ 
can both be identified as Weyl-invariant subrings of $RS$.
We will show that $RG \lt RK$ is surjective in all cases,
so that $\im(RG \to RH \mono RS) = \im(RK \mono RS)$
and hence $(RS)^{WK} \iso RK$ is a subring of $(RS)^{WH} \iso RH$
despite the absence generically of any containment relation between $K$ and $H$.
In fact, one has $WH \leq WK$ as subgroups of $\Aut S$.
To see this,
	\comment{
let $h \in H = G^{\ang{\tau \o c_s}}$ normalize $S$,
so that one has $shs\- = \tau\-(h)$.
Evidently the restrictions to $S$ of the conjugation maps $c_{shs\-}$ and $c_h$ are the same, so
differentiating, we see $\Ad(shs\-)$ and $\Ad(h)$ have the same action on $\fs$,
and selecting $y \in \f h$ such that $\exp y = h$ and differentiating again,
one sees that $\ad\mn\big(\!\Ad(s)y\big)$ and $\ad y$ give the same derivation on $\f s$.
We may form the $\ang\tau$-invariant average
$x = \frac 1 m \sum_{i=1}^m \tau_*^{-i} y = \frac 1 m \sum_{i=1}^m \Ad(s)^i y$ in $\f h$,
and $\ad x$ and $\ad y$ agree on $\f s$
so that $k = \exp x$ lies in $H$ and $c_k$ and $c_h$ yield the same automorphism of $S$. 
	}
recall that the inclusion $N_K(S) \lt N_G(S)$
induces an isomorphism $W_K \lt \pi_0 N_G(S)$~\cite[Prop.~4.4]{goertschesnoshari2014},
and the inclusion $N_H(S) \lt N_G(S)$
induces an injection of reflection groups $W_H \lt \pi_0 N_G(S) \iso W_K$.
Since $\pi_1 G = 0$,
there is a natural identification of $RS$ with the group ring $\Z X$
of the character group $X = \Hom(S,S^1)$,
on which the Weyl groups act as reflection groups, 
so a result of Steinberg~\cite[Thm.~2.2]{steinberg1975pittie}\footnote{\ 
	This theorem explicitly assumes $\pi_1 K = 0$,
	but in the proof of Lemma 3.1,
	specifically the third line from the bottom of p.~175,
	it is pointed out that this result only depends on hypotheses on the root data
	of $K$ also satisfied when $K = \SO(2r+1)$.
}
tells us that $(RS)^{WH} \iso RH$ 
is free over $(RS)^{WK} \iso \im(RG \lt RS)$.
Thus, once we show $RG \lt RK$ is surjective,
we will know $RH$ is free over $\im(RG \to RH)$,
concluding the proof.
}

\rev{
It remains then only to 
show that for diagram automorphisms $\tau$
of simply-connected simple groups $G$,
the map $RG \lt RH$ induced by the inclusion of $H = G^{\ang\tau}$ is surjective
with kernel generated by a regular sequence.
}
\dead{
	Recall~\cite[Prop.~D.40]{fultonharris}
	that every automorphism of a simple Lie group $G$
	is conjugate to one of the form $\tau \o c_t$,
	where $\tau$ is induced by an automorphism of the Dynkin diagram of $G$
	with respect to a maximal torus $T$ 
	and $c_t$ is conjugation by an element of this torus. 
	We assume $\s$ is of this form.
	It is not hard to see that $\defm K = (G^{\ang{\tau}})^0$
	shares a maximal torus $\defm S$ with $H$~\cite[{\SS}4]{goertschesnoshari2016}.
	We will show $\K_K(G/K)$ is of the claimed form.
	Once we have shown this, the result in general will follow
	because \cite[(3.3), p.~176]{snaith1972kunneth}\cite{mcleod1979kunneth}
	\[
		\KS(G/S) 
			\iso 
		\K_{S \x S}  G
			\iso
		R(S \x S) \ox_{R(K \x K)} \K_{K \x K} G
			\iso
		\xt{RG}{RS}{RS} \ox \ewP,
	\]
	and then as in \Cref{rmk:maxrank} one recovers
	$
		\KH(G/H) 
			\iso 
		RH \ox_{RG} RH \ox \ewP,
	$
	as claimed.
}

\dead{
For all $K$ of this form, 
in fact,
$RG \lt RK$ is surjective
with kernel generated by a regular sequence
of length $\rk G - \rk K$
and $RK$ is polynomial, so that the formula 
for $K^*_K(G/K)$ follows from \Cref{thm:TNHZ}.
First note that if $K$ is simply-connected,
then 
$RK$ is a polynomial ring on the fundamental 
representations corresponding to a system of simple roots.
The Dynkin diagram of $K$ with respect to $S$
is the quotient of the Dynkin diagram of $G$ with respect to $T$ 
by the cyclic subgroup generated by a diagram automorphism.
As vertices of each Dynkin diagram correspond to a system of simple roots,
surjectivity of $RG \lt RK$
and the claim about the kernel follow.
This actually covers almost everything, 
as t}\rev{T}he automorphisms 
of Dynkin diagrams
of simple Lie groups comprise the \rev{five} cases
\[
	\big(\SU(2n+1),\SO(2n+1)\mnn\big),\quad
      \big(\SU(2n),\Sp(n)\mnn\big), \quad
\big(\Spin(2n),\Spin(2n-1)\mnn\big),  \quad
(E_6,F_4), \quad
\big(\Spin(8),G_2\big)\mathrlap.
\]
The first four correspond to the classical symmetric spaces of types
AI ($m$ odd), AII, BDI ($q=1$),\footnote{\ 
	Minami calls this case \emph{BDII(a)};
	the ``II'' cuts out the $q=1$ case from the general BDI
	case of $\big(\SO(p+q), \SO(p) \x \SO(q)\big)$ or a cover,
	possibly because the denominator has only one simple factor when $p > q = 1$; 
	the meaning of the``(a)'' is uncertain. 
}
and EIV respectively,
and Minami shows that for each of these pairs $(G,H)$
one has $RG \lt RH$ surjective with kernel generated by a 
regular sequence~\cite[pp.~632, 629--30]{minami1975symmetric}.
For the remaining case, coming from the triality automorphism of $D_4$,
the inclusion $G_2 \longinc \Spin(8)$ factors through a $\Spin(7)$ subgroup.
We can write $R\Spin(8) \to R\Spin(7) \to RG_2$ as a map of polynomial rings
\[
  \Z[\rho,\l^2\rho,\D_+,\D_-] \lt \Z[\s,\l^2\s,\D] \lt \Z[\s,\Ad]\mathrlap,
\]
where
$\rho$ is the standard representation of $\Spin(8)$ factoring through $\SO(8)$
and $\s$ for the analogous representation of $\Spin(7)$ and its restriction to $G_2$,
exterior powers are denoted by $\l^j$,
the adjoint representation of $G_2$ is $\Ad$,
and $\D_\pm$ and $\D$ are respectively the half-spin and spin representations.
It is standard (and shows up in the proof of the BDII(a) case)
that $\rho \lmt \s + 1$ and $\D_\pm \lmt \D$,
so that the kernel of the first map is $(\D_+ - \D_-)$.
It is definitionally true that $\s \lmt \s$ under the second map,
and can be checked~\cite{LiEonline} that $\l^2 \s \lmt \s + \Ad$ and $\D \lmt \s + 1$,
so that the kernel of the second map is $(\D-\s-1)$.
It follows the composite is surjective with kernel generated by the regular sequence
$\D_+ - \D_-, \D_- - \rho$,
concluding the case analysis.
%
%
%
%
%
\epf

\brmk 
Without trading the automorphism
for a diagram automorphism, 
the proof would be significantly more complicated,
\rev{running through the classification~\cite[\selectlanguage{russian}Теорема\!\!\selectlanguage{english}~11]{terzic2001formal} 
of irreducible generalized symmetric pairs $(G,H)$ with $\rk G > \rk H$.
We include the outline of a sample computation.}
\dead{
and we do not know how we would compute the Tor directly.
For an example, compare the formulae in Minami's computation of $\K(G/H)$ 
in the other cases~\cite[e.g., Lem.~5.2]{minami1976symmetric} 
and imagine the relations not in $RH \ox_{RG} \Z$ but in $RH \ox_{RG} RH$. 
It is not at all obvious from that description that things 
should simplify as much as they do.
}
\rev{
The pair $(G,K)$ corresponding to the symmetric space of type FI
is $\big(E_6,\PSp(4)\big)$.
Minami~\cite[p.~632]{minami1975symmetric} observes that $\PSp(4)$ shares a maximal torus $S$ with $F_4$,
which is the fixed point subgroup under the diagram automorphism, 
so that $RE_6 \lt RF_4$ is surjective,
and that both contain a subgroup $Q \iso \SU(2) \ox_{\Z/2} \Sp(3)$
of equal rank.
He then computes the image of $R\PSp(4) \lt RQ$
and the images of three of four generators of $RF_4$ in $RQ$
and is able to conclude that $R\PSp(4) \ox_{RE_6} \Z$
is free abelian of rank $3$.
We do not have the option of reducing modulo $IE_6$
if we want to calculate $\K_K(G/K)$ explicitly,
so we are forced to identify the image of the remaining generator of $RF_4$ in $RQ$,
express this generator in terms of the generators of the image of $R\PSp(4)$,
and show that the one image is free of rank three as a module over the other.
}

\rev{
The ring $RF_4$ is a polynomial ring on the adjoint representation $\Ad$
and the first three exterior powers of an irreducible $27$-dimensional representation $\tau$.
The ring $R\SU(2)$ is $\Z[\sigma]$ for~$\s$ the defining representation,
and the ring $R\Sp(3)$ is polynomial on the first three exterior powers $\rho,\l_2,\l_3$
of the defining representation $\rho$.
  Minami observes~\cite[p.~272]{minami1976symmetric}
that the quotient map embeds $RQ$ in $R\big(\SU(2) \x \Sp(3)\big) = \Z[\s,\rho,\l_2,\l_3]$
as the (non-polynomial) subring $\Z[\l_2,\rho^2,\l_3^2,\rho \l_3,\s^2,\rho\s,\s\l_3]$;
the generators are $\l_2$ and the monomials of degree $2$ in the variables $\rho,\s,\l_3$.
He also observes~\cite[p.~281]{minami1976symmetric} that the image of the injection $R\PSp(4) \lt RQ$
is generated by 
\[\l_2 + \rho\s + 1,\quad
	\s\l_3 + 2\l_2,\quad
	(\rho+\s)^2,\quad
	(\l_3 + \s\l_2 + \rho)^2,\quad
	(\rho+\s)(\l_3 + \s\l_2 + \rho)\mathrlap.
\]
which we will respectively rename $\a,\b,\g^2,\d^2,\g\d$.
More precisely, it is the quotient of the polynomial ring in five indeterminates
by the ideal generated by the single relation $\g^2\d^2 - (\g\d)^2 = 0$.
}

\rev{
	Minami~\cite[p.~280]{minami1976symmetric} computes the image of $\Ad$ in $RQ$ to be $\s\l_3 - \l_2 + \s^2 + \rho^2 - \s\rho - 1$
and those of $\tau$ and $\l^2\tau + \tau $ to respectively be
$\l_2+\s\rho-1$
and
$\rho\l_3 + (\s^2+\rho\s-3)\l_2 +\rho^2$.
He does not determine the image of $\l^3\tau$, 
leaving that task for us. 
Using software to calculate the decompositions of exterior powers
into irreducible representations and then converting nine of these
back into polynomials in $\tau$, $\l^2\tau$, and $\l^3\tau$,
one can work out that $\l^3(\tau+1) = \l^3\tau + \l^2\tau \in RF_4$
is sent under restriction to $Q$
to 
\[\l_3^2 + (\s^2-4)\l_2^2 + \rho^2 \l_2 + (\s\rho^2+\s^3-3\s)\l_3 - \rho^2 + 1\mathrlap.\]
It is not even immediately obvious from these expressions that the image
of $RF_4 \lt RQ$ lies in that of $R\PSp(4) \lt RQ$, though we of course know it from
first principles.
But through dogged persistence in the case of the first three named
generators of $\im(RF_4 \to RQ)$ and Macaulay$2$ computation in the last
case, we are able to express them respectively as
\[
\g\d-\a-\b,\qquad
\a-2,\qquad
\g^2+\b-2\a,\qquad
\theta \ceq \b\g^2+\d^2-2\g\d-\b+2\a-2\a\b+1\mathrlap.
\]
With some rearrangement, we can re-express the image of $RF_4$
as the polynomial ring on 
\[
\a,\qquad
\wt\b \ceq \b+\g^2,\qquad
\vk \ceq \g\d+\g^2,\qquad
\t \mathrlap.
\]
Rewriting $\t$ as far as possible in terms of the new generators,
one finds we can substitute $\wt\t \ceq \d^2 + \g^2(3+2\a+\b)+1$. 
We can then change generators in the image of $R\PSp(4)$
to arrive at $\Z[\a,\wt\b,\vk,\wt\t,\g\d]$.
We want to see that this is free over $\Z[\a,\wt\b,\vk,\wt\t] \iso RF_4$
on the basis $1,\g\d,(\g\d)^2$.
For this we translate the relation $(\g\d)^2 - \g^2\d^2$
into the new generators, using $\g^2 = \vk -\g\d$
and $\d^2 = \wt\t - (\vk - \g\d)(3+2\a+\wt\b + \g\d - \vk)$.
This expands into a monic irreducible cubic in $\g\d$ over $\Z[\a,\wt\b,\vk,\wt\t]$,
finally showing $R\PSp(4)$ is free of rank $3$ over the image of $RE_6$.
}

\rev{One can see from this example the desirability of 
a general principle rather than a case analysis.
To translate this to a statement about $\K_{\PSp(4)}\big(E_6/\PSp(4)\big)$,
recall that $R\Sp(4)$ is generated as a polynomial ring by the first four
exterior powers of the defining representation $\vp$
and the projection $\Sp(4) \lt \PSp(4)$ embeds $R\PSp(4)$
as the subring of $\Z[\vp,\l^2\vp,\l^3\vp,\l^4\vp]$ generated by
$\l^2\vp, \l^4\vp, \vp^2, (\l^3\vp)^2, \vp\l^3\vp$,
whose images are respectively the elements $\a,\b,\g^2,\d^2,\g\d$ from earlier.
Write $z$ and $w$ for the images under Hodgkin's~\cite{hodgkin1967lie}
map $QRG \lt P\K G$ of the differences of fundamental
representations $\pi_1 - \pi_2$ and $\l^2\pi_1 - \l^2\pi_2$
corresponding to the simple roots identified by the $\Z/2$ symmetry
of the Dynkin diagram of $E_6$.
}

\rev{
Then finally we can write
$
\K_{\PSp(4)}\big(E_6/\PSp(4)\big)
$ 
as the quotient of 
\[\Z[\l^2\vp,\  \l^4\vp,\  \vp^2,\  (\l^3\vp)^2,\  \vp\l^3\vp,\ 
\l^2\psi,\  \l^4\psi,\  \psi^2,\  (\l^3\psi)^2,\  \psi\l^3\psi] \ox \ext[z,w]\]
by the ideal
\begin{multline*}
\big(
\l^2\vp-\l^2\psi,\ \  \l^4\vp +\vp^2-\l^4\psi-\psi^2,\ \ 
\vp^2+(\l^3\vp)^2-\psi^2-(\l^3\psi)^2,\\
\vp\l^3\vp + \vp^2(3+2\l^2\vp+\l^4\vp)-
\psi\l^3\psi - \psi^2(3+2\l^2\psi+\l^4\psi)
\big)
\mathrlap.
\end{multline*}
}
\end{remark}

\begin{example}
	There is one remaining class of pairs $(G,H)$ 
	for which \eqfity of the isotropy action is characterized in a computationally
	tractable way.
	If $H$ is a circle, the \iact is equivariantly formal if and only if
	either $\pi_1 H \lt \pi_1 G$ is injective or (these options are exclusive)
	there is an element $g \in G$ such that conjugation by $g$ 
	induces the nontrivial automorphism $h \lmt h\-$ of $H$.
	
	In the former case, \Cref{thm:TNHZ} 
	applies after inverting the order $\ell$ 
	of the center of the commutator group $G'$,
	since the composition $H \inc G \epi G/G' = G\ab$ is then at most $\ell$-to-one,
	so $R(G\ab)[1/\ell] \lt R(H;\Q)[1/\ell]$ and hence 
	$R(G)[1/\ell]\lt R(H;\Q)[1/\ell]$ are surjective.
	(Consider the diagonal circle group in $\U(2)$ to see this inversion is necessary.)

	In the latter case, it is possible that
	$RH$ not be flat over the image of $RG \lt RH$,
	so the current proof cannot be adapted to handle this case.
	For example, \cite[Ex.~7.18]{carlsonfok2018}, 
	for $G = \SU(4)$ and $H = \diag(z,z\-,z^2,z^{-2})$, 
	the Laurent polynomial ring $RH \iso \Z[t,t\-]$ 
	is not flat over the image of $RG \lt RH$.
	%
	%
\end{example}

\brmk\label{rmk:Abhyankar}
	The strength of the hypothesis that the image $R$ of $RG \lt RH$
	be a complete intersection ring in \Cref{thm:main}
	has been an irritating limitation.
	The hypothesis would not be as seemingly overpowered if it were redundant,
	which seems \emph{a priori} reasonable.
	That is, it is possible the following holds:

\medskip

\emph{For any 
	surjection 
	$\varphi\: A \longepi B$
	of polynomial rings, respectively
	in $m \geq n$ indeterminates 
	over a commutative base ring $\ring$, 
	one can choose an algebraically independent 
	set $x_1,\ldots,x_n,y_{n+1},\ldots,y_m$
	of $\ring$-algebra generators for $A$
	such that $\varphi$ sends $y_j \lmt 0$ 
	and restricts to an isomorphism $\ring[x_1,\ldots,x_n]\isoto B$.
}

\medskip

This claim
is easily verified for graded maps of graded rings over $\ring = \Q$,
leading to the analogous result for maps $H^*(BG;\Q) \lt H^*(BH;\Q)$,
but is a hard open problem for ungraded maps even over $\ring = \C$.
In geometric language, 
the special case $m = n+1$
is the \emph{Abhyankar--Sathaye embedding conjecture}%
~\cite{abhyankarmoh1975,%
	sathaye1976linear,%
	russellsathaye,%
	popov2015around,%
	MO:Wendt},
which states that any embedding 
$\mathbb{A}^{\mn n}_\C \longmono \mathbb{A}_\C^{\mn n+1}$
is taken to the standard embedding 
by some automorphism of $\mathbb{A}_\C^{\mn n+1}$.
This is known at present for $n = 1$ and several other special cases, 
and is closely related to the determination of the algebraic automorphism group
$\Aut \mathbb{A}_\C^{\mn m}$, which is still incomplete for $m \geq 3$.

The author knows of no instance where the claim does not hold
for surjections $RG \lt RH$ between representation rings 
of compact, connected Lie groups, but also does not know any 
specific characteristics of these rings that make the claim 
any easier. 
In the case $G/H$ is an odd-dimensional sphere,
the claim can be shown case-by-case using the classification,
and the fact $RG$ and $RH$ are $\l$-rings is helpful in these cases,
but there do not seem to be known properties of 
$\l$-ring homomorphisms that would generalize this proof.

One can also attempt to prove the result using a classification 
of all pairs $(G,H)$ where $RG \lt RH$ is surjective
and $\pi_1 G$ is torsion-free.
There are a few useful initial reductions before
this approach runs into difficulties.
With an additional condition regarding
the embedding of the center of $H$ in $G$,
one can reduce to the case of $G$ semisimple,
and hence, by the $\pi_1$ condition, a product of simply-connected
simple groups.
If one assumes $H$ as well is a product of simple factors,
the surjectivity of $RG \lt RH$
reduces to the surjectivity of $\Tensor RG_j \lt RH_\ell$
where $H_\ell$ is a simple factor of $H$ contained in a product 
$\prod G_j$ of simple factors of $G$.
But the image of this map is the subring generated 
by the images of the restrictions $RG_j \lt RH$,
none of which individually is surjective,
so to obtain a surjective such map,
one need only pick a collection of 
representations of $H$ generating $RH$
and for each such representation $V_j$ 
pick some simple group $G_j \leq H$
to which $V_j$ extends.
Thus the classification 
of surjections of representation rings $RG \lt RH$
for $(G,H)$ compact connected
is essentially equivalent
to the problem of describing all maps
of representation rings of simple groups.
\ermk

\section{The map of spectral sequences}\label{sec:mapkunneth}
This section collects definitions, constructs Hodgkin's spectral sequence
and the promised map, and provides a cautionary counterexample and historical context.

	To conduct a comparison of {\KSS}s, 
	we need to analyze their construction.
	Whereas the original
	construction of the \EMSS of that square starts from the isomorphism
	$H_G^*(X \x Y;\Z) \iso \Tor_{C^*(BG)}\mn\big(C^*(X_G),C^*(Y_G)\mnn\big)$,
	requiring a notion of ``differential Tor'' of a differential module over a \DGA,
	and proceeds through an algebraic K{\"u}nneth spectral sequence,
	an analogous construction in equivariant K-theory
	is unavailable due to the general lack of a cochain-level model.
	A later method of constructing the \EMSS
	developed by Larry Smith~\cite{smith1970emss},
	based on ideas of himself, Hodgkin, and Rector,
	navigates around this difficulty by proceeding first at the space level, 
	constructing a ``geometric resolution'' 
	in an appropriate category of topological spaces, and obtaining the algebraic resolution through functoriality.
	Hodgkin showed these ideas apply to a broad class of cohomology theories~\cite[{\SS}1--6]{hodgkin1975kunneth}.

\begin{construction}
	Let $(\ms C,\rev{\unit}\dead{*},\^)$ be a pointed monoidal category
	with finite products and colimits,
	tensored over pointed CW complexes
	in such a way that there is an isomorphism 
	$\rev{\unit}\dead{*} \^ S^0 \iso \rev{\unit}\dead{*}$.\footnote{\
		In a typical application objects of $\ms C$ will be topological spaces
		with some extra structure and the tensor operation will be smash product.
	}
	We equip $\ms C$ with a notion of a \defd{homotopy}
	of a morphism $f\: X \lt Y$,
	given as a map $X \^ I_+ \lt Y$,
	and a notion of a \defd{cofiber},
	given as the pushout of $Y \from X \to X \^ I$.\footnote{\ 
		This notion of homotopy is clearly related 
		to the notion of left homotopy with respect to
		a functorial cylinder object $X \^ I$ for
		a model category structure in the intended applications,
		but we have made no attempt to find an exact translation
		of the minimal hypotheses in terms of the constellation
		of notions surrounding cofibration categories,
		Waldhausen categories, etc.
		}
	Let $\defm{h^*}$ be a reduced multiplicative cohomology theory on $\ms C$,
	meaning a contravariant (commutative graded algebra)--valued homotopy functor
	$h^* = \Direct_{n \in \Z}  h^n$ 
	taking a cofiber sequence
	to an exact sequence of groups
	and equipped with an element 
	$\defm{\es} \in  h^1 (\rev{\unit}\dead{*} \^ S^1)$
	such that the map
	$ h^* X \isoto  h^{*+1} (X \^ S^1)$ 
	given by the cross product with $\es$
	agrees with the suspension isomorphism
	given by the cofiber sequence $X \to X \^ I \to X \^ S^1$.
	

	An object $Z \in \ms C$ is called a \defd{K\"unneth object} 
	for ${h}^*$
	if
	$ h^* Z$ is a finitely-generated projective module over
	the coefficient ring 
	$\defm{h^*} \ceq  h^*(\rev{\unit}\dead{*} \^ S^0)$ 
	and the cross product
	${ h^* Z} \ox_{h^*}{ h^*Y} \lt  h^*(Z \wedge Y)$
	is an isomorphism natural in $Y$.
	A \defd{geometric resolution} of $X$ with respect to $ h^*$ is a
	sequence 
	\[
	X 
	\eqc X_0 \lt Z_0 \lt X_1 \lt Z_1 \lt X_2 \lt Z_2 \lt \cdots
	\]
	in $\ms C$ 
	such that each $Z_\ell$ is a K\"{u}nneth object, 
	each $h^* Z_\ell \lt h^* X_\ell$ is surjective,
	and each sequence
	\quation{\label{eq:cofiber}
		X_\ell \os{j_\ell}\lt Z_\ell \os{\kk_\ell}\lt X_{\ell+1},
	}
	is a cofiber sequence. 
	Tensoring \eqref{eq:cofiber} with $Y$ yields further cofiber sequences,
	yielding long exact sequences on applying $ h^*$,
	and adding them gives an exact couple
	\[
	\xymatrix@R=2.25em@C=-2em{
		\Direct  h^*(X_\ell \wedge Y)	\ar[rr]^(.475)\d
		&& \Direct  h^*(X_{\ell+1} \wedge Y) 	\ar[ldd]^j\\ \\
		&	\Direct  h^*(Z_\ell \wedge Y).\ar[luu]^k
	}
	\]
	The \defd{K\"{u}nneth spectral sequence} for $\thy$ is 
	the cohomological left-half-plane spectral sequence associated to this exact couple~\cite[p.~283]{snaith1971massey}.
	The composition $jk\:  h^*(Z_\ell \wedge Y) \lt  h^*(Z_{\ell-1} \wedge Y)$
	can be identified with 
	$ h^* Z_\ell \ox_{h^*}  h^*Y \lt  h^* Z_{\ell-1} \ox_{h^*}  h^* Y$,
	by the K{\"u}nneth property of the $Z_\ell$,
	and since we took $\kk_\ell^*$ surjective,
	we can piece together the cohomology of the sequences (\ref{eq:cofiber})
	into an $h^*$-module resolution of $h^*X$:
	\[
	\cdots \to  h^* Z_2  \to  h^* Z_1 \to  h^* Z_0 \to h^* X \to 0.
	\]
	It follows the $jk$-cohomology of $E_1 = \Direct  h^*(Z_\ell \wedge Y)$ 
	is $E_2 = \Tor_{h^*}( h^* X,  h^* Y)$. 
	The spectral sequence strongly converges~\cite[Thm.~5.1]{hodgkin1975kunneth},
	but at this level of generality it is impossible to say 
	whether it is to $ h^*(X \wedge Y)$ as hoped.\footnote{\ 
		This construction is enough to obtain the bigraded group structure
		on the spectral sequence and discuss convergence issues,
		but the product structure is only provided by the definition
		in terms of a Cartan--Eilenberg $H(p,q)$-system, 
		which we choose not to discuss here.}
\end{construction}

\begin{discussionsilent}
	Smith~\cite{smith1970emss} showed this technique recovers the \EMSS
	when applied to the functor $\thy(X \to B) \ceq \wt H^* X$ 
	on the category $\smash{\Top_{/B}}$ of spaces over a fixed space $B$,
	and Hodgkin applied it to $\thy = \wt K^*_G$ on 
	the category $G\mbox{-}\Top$ of $G$-spaces.
	To construct geometric resolutions 
	Hodgkin needs to prove
	there are ``enough'' K\"unneth spaces for $\KG$,
	and finds products of suspensions of 
	complex Grassmannians $\Gr(\ell, V)$ of
	complex $G$-representations $V$ will do%
	~\cite[Prop.~7.1]{hodgkin1975kunneth}.
	This is itself a substantial technical result.

	To see there is a homomorphism of spectral sequences,
	note that these $G$-spaces $\Gr(\ell,V)$
	are also K\"unneth spaces for $H_G^{**}(-;\Q)$,
	then we may simultaneously take the same geometric resolution for
	both theories,
	and apply the equivariant Chern character 
	to these resolutions
	to induce a homomorphism of K\"unneth spectral sequences. 
	But it follows from \eqref{eq:extension} that any K\"unneth space for $\KG$
	is also one for $H_G^{**}(-;\Q)$.\footnote{\ 
		Alternately, one can prove directly that given a $G$-representation $V$,
		the $G$ action on $Y = \Gr(\ell,V)$ is equivariantly formal 
		(much easier than what Hodgkin had to do) 
		which implies $H_G^{**} Y$ is free over $H_G^{**}$,
		so that the 
		\KSS associated to the right square in \eqref{eq:Gfibersquare}
		collapses to a K\"unneth isomorphism.
		For this, let $n = \dim_\C V$
		and $H = \U(\ell) \x \U(n-\ell)$.
		Then $\Gr(\ell,V) = \U(n) / H$ with the left multiplication action of $G$,
		and the action of $G$ on $\U(n)/H$ is equivariantly formal if and only if the map
		$\chi\: \U(n)/H \to \big(\U(n)/H\big)_G$ induces a surjection.
		
		First, this is so if $G = \U(n)$,
		for then $\chi\: G/H \lt BH$ is 
		the classifying map of the principal $H$-bundle $G \to G/H$,
		and $\chi^*$ is surjective since $\rk \U(n) = \rk H = n$.
		Second, it is true if $G < \U(n)$
		since $\U(n)/H \lt E\U(n) \ox_{\U(n)} \U(n)/H$
		factors through $E\U(n) \ox_G \U(n)/H$.
		Finally, in general, let $\G < \U(n)$ be the image of $G$
		and $N \normal G$ the ineffective kernel.
		Then there is a free right action of $\G = G/N$ on $BN \simeq EG/N$,
		so 
		\[
		EG \ox_G \U(n)/H \homeo BN \ox_\G \U(n)/H = \big(BN \x \U(n)/H\big)_\G\mathrlap,
		\]
		but $\U(n)/H \to \big(\U(n)/H\big)_\G$ factors through
		$\U(n)/H \to \big(BN \x \U(n)/H\big)_\G$.
	}
\end{discussionsilent}

%

\SPECMAP*

\brmk
The comparison of the \HKSS with more tractable sequences
has precedent.
Hodgkin's procedure applied to the theory $(X \mn\mn\to\mn\mn BG) \longmapsto \K  X$ 
on $\Top_{/BG}$
yields a spectral sequence starting at $E_2 = \Tor_{\K BG}(\K  X, \K  Y)$
and converging to $\K( X \x_{BG}  Y)$ when $\pi_1 G$ is free abelian,
and he uses it in an essential way in his convergence proof (which in its original
version only shows the $\KG$ sequence converges up to completion)%
~\cite[{\SS}8]{hodgkin1975kunneth}.
Snaith~\cite[{\SS}5]{snaith1971massey}, without considering a map of {\KSS}s,
examines the {\AHSS}s $H^*(X_G;\Z) \implies \K X_G$ 
and $H^*(BG;\Z) \implies \K BG$
to show that if $H^*(BG;\Z)$ and $H^*(X_G;\Z)$ are torsion-free
and the \EMSS $\Tor_{H^*(BG;\Z)}\big(\Z,H^*(X_G;\Z)\mnn\big) \implies H^*(X;\Z)$
collapses, then so also does the \HKSS
$\Tor_{RG}(\Z,\KG X) \implies \K X$,
which we will cite as \Cref{thm:Snaith-inj}.
\ermk

\brmk\label{rmk:convtorsion}
The demand on $\pi_1 G$ arises from the typical convergence 
of the \HKSS
to a target which is distinct from $\KG(X \x Y)$ 
in the contrary case. 
Hodgkin demonstrates this failure with two 
examples~\cite[p.~68]{hodgkin1975kunneth}.
The first is $G = X = Y = \Z/2$, 
which has $\smash{E_2 = \Tor_{R(\Z/2)}(\Z,\Z) = H_*(B\Z/2)}$
equal to $0$ in even positive degrees and $\Z/2$ in odd degrees,
and collapses for lacunary reasons.
On the other hand, 
$\K_G(G \x G) = K^0(\Z/2) = \Z^2$. 
The second is $G = X = Y = \SO(3)$;
here $E_2 = \Tor_{R\SO(3)}\mnn(\Z,\Z)$ is an exterior
algebra on one generator since $R\SO(3)$ is a polynomial ring on one generator,
but $\K_G(G \x G) = \K\SO(3)$ contains 2-torsion.
\ermk

\bex\label{ex:SO(3)}
One could be forgiven for hoping from Hodgkin's examples of bad convergence 
(and the wording he uses) that the problem for a connected $G$ 
might lie solely with torsion: 
after all, the \KSS in Borel cohomology requires $G$ to be connected too,
and cohomology with $\Q$ coefficients is not very sensitive to torsion in $\pi_1$ 
(e.g., replacing $(G,H)$ with a connected finite cover does not 
affect \eqfity of the \iact~\cite[Thm.~1.2]{carlson2018eqftorus}).
This hope is in vain. 

Write $G = \SO(3)$ and $H = T = \SO(2)$,
so that the isotropy action of $H$ on $X = Y = G/H = S^2$
is the standard rotation of a globe.
One would like to use the \KSS to compute $\KT(G/T)$
starting from 
\[
E_2 = \Tor_{RG}\mn\big(\KG(G/T),\KG(G/T)\mnn\big) 
\iso
\Tor_{R\SO(3)}\!\big(\mn R\SO(2),R\SO(2)\mnn\big).
\]
Writing $\defm t \: \SO(2) \isoto \U(1)$
for the standard representation,
we have $R\SO(2) = \Z[t,t\-]$
and $R\SO(3) \iso \big(\mnn R\SO(2)\mnn\big)\mn{}^{W\SO(3)} = \Z[t+t\-]$,
so particularly $R\SO(2)$ is a free module of rank $2$ over $R\SO(3)$,
say on $1$ and $t$,
and $\Tor = \Tor^0$ is the tensor product $R\SO(2) \ox_{R\SO(3)} R\SO(2)$,
free of rank two over $R\SO(2)$.
Thus the \KSS must collapse.
We have in general, for $G$ compact and connected 
and $T$ its maximal torus, a sequence of maps~\cite[1.4]{mcleod1979kunneth}
\[
\xt{RG}{RT}{RT} 
\os{\defm\l}\lt
\KT(G/T)
\os{\defm\iota}\lt
\KT W,
\]
where $\l$ is the edge map
of the spectral sequence---%
which is an isomorphism if the collapsed sequence converges to $\KT(G/T)$---%
and $\iota$ is restriction to the fixed point set
$(G/T)^T = N_G(T)/T = W$ of the action.
It is always the case that $\iota$ is an injection~\cite[Thm.~1.6]{mcleod1979kunneth},
so to show $\l$ is not an isomorphism in the case at hand it will be enough 
to see the images of $\l$ and $\l \o \iota$ differ.

We can calculate $\iota$ directly 
from the Mayer--Vietoris sequence corresponding to the cover of $S^2$ 
by invariant hemispheres
meeting in an equatorial $S^1$-orbit,
since these each contract equivariantly to a fixed point.
Since $K^1_{S^1}(S^1) = 0 = K^1_{S^1}(*)$, 
the sequence reduces to 
\[
0 \to K^0_{S^1}(S^2) \os{\iota}\lt RS^1 \x RS^1 \lt \Z \to 0,
\]
where both maps $RS^1 \lt \Z$ are the augmentation.
The kernel of their difference, $\im \iota$, is thus given by
pairs of virtual representations of equal dimension.
Since $t^n$ for $n \in \Z$ all are one-dimensional,
$\im \iota$ admits as a $\Z$-basis 
the pairs $(t^m,t^n) \in RS^1 \x RS^1$ 
for $n,m \in \Z$
and as ring generators the four pairs $(t^{\pm 1},1)$ and $(1,t^{\pm 1})$.
On the other hand, $\l$ takes
%
\eqn{
	1 \ox t &\ \lmt\  \big[\xt{\SO(2)}{\SO(3) } \C \to \xt{\SO(2)}{\SO(3)}{\{0\}}\big],\\
	t \ox 1 &\ \lmt\  [S^2 \x \C \to S^2],
}
where $\SO(2)$ acts on the left of $\SO(3)$ in the first bundle 
and diagonally in the second.\footnote{\ 
	These expressions are less mysterious than they may seem. 
	They come from expanding the definition of the action of $RH \x RH$ on $\K_{H \x H} G$
	and the isomorphism $\K_{H \x H} G \isoto \KH(G/H)$.
}
Orienting $S^2$ and restricting these bundles to their representations at the poles,
we get respectively $(t,t)$ and $(t,t\-)$,
showing $\im \iota$ is of rank two over $\im(\iota \o \l)$.

As these groups are all torsion-free and 
$\SO(3)$ is the simplest group that is connected with fundamental group 
containing torsion,
we see there is no hope of saving convergence by extending coefficients.
This example also shows \Cref{thm:main} fails, with arbitrary coefficients,
without the assumption $\pi_1 G$ be torsion-free.
\eex

\brmk\label{rmk:convergence}
The convergence proof 
for the \HKSS
follows a long series of reductions which we present
for that portion of posterity with the stamina;\num{742}
the proof combines the work of several authors 
in non-chronological order, 
and it seemed worthwhile that the logical structure
of the proof should be outlined somewhere.
It is also useful to see where the insistence on 
torsion-free $\pi_1 G$ comes from, 
and further, 
the proof of the theorem motivates the counterexample~%
\ref{ex:SO(3)}.

First, one sees the \HKSS converges to a colimit
$F(X,Y) \ceq \colim \wt K^*_G\big(({\susp^p X_+}/{X_p}) \, \wedge\, Y_+\big)$,
independent of the resolution chosen~\cite[Thm.~5.1 \& p.~42]{hodgkin1975kunneth},
and this colimit sits in an exact triangle 
\[
F(X,Y) \to \KG(X\x Y) \to \G(X,Y) \to F(X,Y)\mathrlap,
\]
so we will be done when we show the error term
$\smash{\defm\G(X,Y) = \colim \wt K^*_G(X_p \wedge Y_+)}$ vanishes identically.
For this note~\cite[p. 60\emph{f}.]{hodgkin1975kunneth} that for fixed $X$,
essentially because the same facts hold for $\KG$,\num{751}
the functor $\defm{\G^X} = \G(X,-)$ 
is a cohomology theory in $Y$, continuous in the sense
that for $A \sub X$ closed, $\G^X A$ is the colimit of $\G^X B$
over closed $B \supseteq A$.
Writing $\pi\: Y \lt Y/G$ for the quotient maps,
we~\cite[Prop.~5.3]{segal1968equivariant} 
have a Segal spectral sequence 
\[
E_2 = H^*(Y/G;\G^X \pi\-) \implies \G^X Y\mathrlap,
\]
where $\G^X \pi\-$ denotes the sheaf with stalks $\G^X \pi\-(yG)$,
and it follows that to show $\G = 0$
we need only show $\G^X(G/H) = 0$ for all $X$ and all closed $H < G$.
Applying the same argument to $\G^{G/H}$ in turn,
to conclude $\G = 0$ we only need to see $\G(G/H,G/K) = 0$ 
for all closed subgroups $H,K < G$.

As a preliminary step, it will be important to show $\G(G, -)$ is $0$.
Since $\pi_1 G = 0$, 
the short exact sequence $1 \to G' \to G \to G\ab \to 1$
induces an unnatural decomposition
$RG \iso RG' \ox RG\ab$,
where the first factor is a polynomial ring and the second a 
Laurent polynomial ring~\cite[Prop.~11.1]{hodgkin1975kunneth}.
The global dimension of such a ring is $1 + \rk G$%
~\cite[Cor.~11.2]{hodgkin1975kunneth},
so projective resolutions are finite.
It follows~\cite[Lem.~8.5 \emph{et seq.}]{hodgkin1975kunneth} that $\G(G,Y)$ is a quotient of $\wt K^*_G(X_{2+\rk G} \x Y,Y)$,
hence finitely generated and discrete in the $IG$-adic topology%
~\cite[p.~67]{hodgkin1975kunneth}.
The homotopy quotient process $X \lmt X_G$
takes a geometric $\KG$-resolution of $X$ for $\KG$
to a geometric resolution of $X_G$ for 
the theory $(X \mn\to\mn BG) \longmapsto \K X $ on $\Top_{/BG}$,
since $\K X_G \iso \KG(X)\compl$,
and for any $\ol Y \in \Top_{/BG}$
we have an exact triangle
\[
\ol F(X_G,\ol Y) \to K^*(X_G \x_{BG} \ol Y) \to \ol \G(X_G,\ol Y)
\mathrlap,
\]
such that if $\ol Y = Y_G$ for some $Y \in G\mbox{-}\Top$,
then $\ol\G(X_G,Y_G) \iso \G(X,Y) \ox_{RG} {\wh RG}$~\cite[p.~64]{hodgkin1975kunneth}.
Because $\G(G,Y)$ is finitely generated and discrete,
we have $\ol \G(EG,Y_G) \iso \G(G,Y)\compl \iso \G(G,Y)$,
so it is more than enough to show the theory $\ol \G(EG,-)$ is zero.
For this we can apply Dold's theorem~\cite{dold1970chern} that 
a natural transformation $e^* \lt f^*$
of additive cohomology theories on $\Top_{/BG}$
is an isomorphism if every map $* \lt BG$ induces an isomorphism.
Since $BG$ is path-connected, 
the inclusion of any point will do,
so it is enough to check that
\[
\ol\G(EG,*) = \ol\G(EG,EG) = \G(G,G)\compl = \G(G,G) = 0
\mathrlap.
\]
To see this,
note that the Laurent-tensor-polynomial structure on $RG$
and the fact $\K G$ is an exterior algebra on the primitives $P\K G$
lets us use the Koszul algebra $RG \ox \K G$ as an $RG$-resolution on $G$,
since Hodgkin's natural transformation $\b\: RG \lt K^1 G$,
defined by viewing a representation $\rho$ as a continuous map 
$G \to \U(n) \inc \U$,\footnote{\ Geometrically,
	this is represented by the equivariant 
	$\C^n$ bundle on the unreduced suspension $SG = CG \union_G CG$  
	obtained by clutching two $G$-trivial bundles on the cone $CG$ 
	via $\rho$ along the common copy of $G$.}
induces an isomorphism $QRG \isoto P\K G$~\cite{hodgkin1967lie}.
Then we can calculate 
\[
E_2 = \Tor_{RG}(\Z,\Z) \iso \K G \iso \KG(G \x G)
\mathrlap,
\]
so that $\G(G,G) = 0$~\cite[p.~85]{hodgkin1975kunneth}.\footnote{\ 
	This is a simplification; one needs to check that the map 
	$\Tor^1_{RG}(\Z,\Z) \to F(G,G) \to K^1_G(G \x G)$ is actually 
	the abstract isomorphism cited.}

For $G = T$ a torus, Snaith~\cite[{\SS}2]{snaith1972kunneth} 
is able to prove $\G(T/H_1,T/H_2) = 0$
for all closed subgroups $H_1,H_2$ of $T$,
using primarily the fact that $\G(T,-) = 0$ since $\pi_1 T$ is free abelian
and homological algebra with respect to specific resolutions that work due to
the restricted nature of subgroups of a torus.
Appending subscripts to label the acting group, this shows
$\defm{\l}_T\: F_T(-,-) \lt \KT(- \x -)$
is a natural isomorphism.
Now Pittie~\cite{pittie1972homogeneous,steinberg1975pittie} 
showed that if $\pi_1 G$ is free abelian 
and $T$ is the maximal torus of $G$,
then $RT$ is a free module over $RG$.
It follows in the \KSS for $G$ acting on $X = Y = G/T$,
we have collapse at $E_2 = \Tor_{RG}(RT,RT) = RT \ox_{RG} RT$.
But McLeod~\cite{mcleod1979kunneth} showed later
in posthumously published work, 
via mainly Lie-theoretic methods, that the map
\[
\l_G\: RT \ox_{RG} RT \to \KG(G/T \,\mn\x\,\mn G/T) \iso\KT(G/T)
\] 
from \Cref{ex:SO(3)} 
is an isomorphism if $\pi_1 G$ is free abelian.
It follows from repeated application of the tricks
$Z \ox_T G \homeo Z \x G/T$
and $\KG(Z \x G/T) \iso \KT Z$
for $Z$ a $G$-space, 
using Pittie's and McLeod's theorems,
that a K\"unneth space for $\KG$ is also one for $\KT$. 

Consider now a $\KG$-resolution $G/T \to Z_0 \to Y_1 \to \cdots$ of $G/T$.
Now~\cite[{\SS}3]{snaith1972kunneth}
since $\KT(G/T) \iso RT \ox_{RG} \KG(G/T)$ by McLeod's theorem
and the $Z_n$ are K\"unneth spaces for $T$ as well,
by induction $G/T \to Z_0 \to Y_1 \to \cdots$ is also a $\KT$-resolution.
For a $G$-space $X$, then,
restriction of actions induces a commutative square
\[
\xymatrix@C=1.5em@R=5.5em{
	F_G(X,G/T) \ar[r]^(.46){\l_G} \ar[d] 	& \KG(X \x G/T)\phantom{.} \ar[d]\\
	F_T(X,G/T) \ar[r]^(.46){\l_T}			& \KT(X \x G/T).
}
\]
We have seen $\l_T$ is an isomorphism, 
and the inclusion $T \inc G$ 
induces ``index'' retractions of the vertical maps%
~\cite[Prop.~3.8]{segal1968equivariant}\cite[Prop.~4.9]{atiyah1968bott},
implying $\l_G$ is an isomorphism as well in this case.
But $\KG(X \x G/T) = \KT X$
and since $RT$ is free over $RG$, 
we have a collapse 
$\Ei = E_2 = \Tor_{RG}(\KG X,RT) = \KG X \ox_{RG} RT = \gr F_G(X,G/T)$,
so we have proven that $\KT X \iso \KG X \ox_{RG} RT$.

Now consider a resolution $X \to Z_0 \to X_1 \to \cdots$ 
of any $G$-space $X$.
Since $\KT X = \KG X \ox_{RG} RT$ and the $Z_n$ are 
also K\"unneth spaces for $T$,
by induction this is also a $T$-resolution, 
so there is an induced commutative square
\[
\xymatrix@C=1.5em@R=4.5em{
	F_G(X,Y) \ar[r]^(.46){\l_G} \ar[d] 	& \KG(X \x Y)\phantom{.} \ar[d]\\
	F_T(X,Y) \ar[r]^(.46){\l_T}			& \KT(X \x Y).
}
\]
But $\l_T$ is an isomorphism and the index induces a retraction, 
so $\l_G$ is an isomorphism as well.

\ermk

\section{The theorem of Shiga and Takahashi}\label{sec:ST}

In our discussion of hypotheses in the introduction,
we adverted to a modification, 
announced without proof 
in the prequel~\cite[Rmk.~3.12]{carlsonfok2018}, 
of a theorem
of Shiga--Takahashi~\cite{shigatakahashi1995}.
We are not able to prove anything like it in the K-theoretic 
case owing to the lack of a natural grading on the representation
ring, but
we prove the cohomological version here. 

\rev{\emph{
In this section, all unadorned tensor products are taken over $\Q$.
}}

\begin{definition}
\dead{In this section, all unadorned tensor products are taken over $\Q$.}
If \GK is a \ccpair such that
the isotropy action of $K$ on $G/K$ 
is equivariantly formal
  we will for brevity call the pair $(G,K)$ \defd{isotropy-formal}.\footnote{\ 
	For consonance with the previous sections
	one might well prefer the isotropy subgroup to be $H$,
	but to harmonize with previous work and also to not use ``$H$'' both
	for this group and cohomology, we prefer $K$ in this section.
  }
If $G/K$ is formal
\rev{in the sense that its minimal model is quasi-isomorphic to its cohomology
ring},\num{Def. 5.1}
we will call the pair~\GK \defd{formal} as well.
Note that $BK = EG/K$ admits a right $N_G(K)$-action given by
$eK\.n = enK$
and $BG = EG/G$ admits a trivial $N_G(K)$-action
making the quotient map $BK \lt BG$ equivariant.
It follows the ring map $H_G^* \lt H_K^*$
is also equivariant, so its image lies in the invariant subring 
$\defm{\HKN} \ceq (H_K^*)^{N_G(K)}$. 
This trivial observation 
explains our apparent fixation on invariant subrings.
Because the action of the identity component of $N_G(K)$ on $BK$ is homotopically
trivial, the group effectively acting on $H_K^*$ is the component group
$\defm N \ceq \pi_0 N_G(K)$.
\end{definition}

The original Shiga--Takahashi theorem is the following:

\begin{theorem}[{\cite[Thm.~A, Prop.~4.1]{shiga1996equivariant}%
		\cite[Thm.~2.2]{shigatakahashi1995}}]
	If \GK is a formal pair and $H_G^* \lt \HKN$ is surjective, 
	then \GK is \isotf.
	If \GS is a formal pair 
	where $S$ is a torus containing regular elements of $G$,
	then \GS is \isotf if and only if $H_G^* \lt \HSN$ is surjective.
\end{theorem}

The published proof of the forward implication 
is somewhat obscure and deserves amplification. 
It turns out the proof of the backward direction can 
be finessed into an equivalence.
In fact, we will show the following, and then derive further consequences:

\begin{restatable}{theorem}{STP}\label{thm:ST+}
	Let $(G,K)$ be a formal pair.
	Then $(G,K)$ is \isotf if and only if $\rho^*: H_G^* \lt \HKN$ is surjective.
\end{restatable}

To get there, we will cite a few lemmas about equivariant formality of
an isotropy action, formality, and classical invariant theory.

\begin{lemma}[{\cite[Prop.~3.1, p.~81]{goertsches2012isotropy}}]%
	\label{thm:fpdim}
	A \ccpair $\GK$ is \isotf
	if and only if $\dim_\Q H^*(G/K) = |N| \. 2^{\rk G - \rk K}$.
\end{lemma}

\begin{theorem}[\textup{\cite[Thm.~1.1]{carlson2018eqftorus}}]\label{thm:SGS}
	If \GK is a \ccpair and $S$ any maximal torus of~$K$, 
	then $(G,K)$ is isotropy-formal if and only if $(G,S)$ is.
\end{theorem}

\blem[{\cite[Rmk., p.~212]{onishchik}
	}]%
\label{thm:formaltorusreplacement}
	If $(G,K)$ is a \ccpair and $S$ a maximal torus of $K$,
	then $G/K$ is formal if and only if $G/S$ is.
\elem

\bthm[{\cite[Thm.~A]{carlsonfok2018}}]\label{thm:isotfformal}
Let $(G,K)$ be an \isotf pair.
Then $G/K$ is formal.
\ethm

Formality itself admits the following formulation for homogeneous spaces,
where we write  $A^{\defm{\geq 1}}$ for the augmentation ideal in a connected graded algebra $A$ and given a map $f\: A \lt B$ of connected graded algebras,
set $\defm{B \quot A} \ceq B / f(A^{\geq 1})B$.

\bthm[{\cite[Thm.~IV, p.~463]{GHVIII}\cite[Thm.~2, p.~211]{onishchik}}]\label{thm:formalCDGA}
Let $(G,K)$ be a \ccpair. Then $G/K$ is formal if and only if there is an isomorphism
\[
H^*(G/K) \iso H_K^* \quot H_G^* \,\ox \ewP,
\]
where the exterior factor $\ewP$ is generated 
by a vector subspace $\wP$ of dimension $\rk G - \rk K$ 
and the ideal generated by $\im(H_G^{\geq 1} \to H_K^*)$
can be generated by a regular sequence of elements of this image.
\ethm

Very briefly, the relevance of the rings in \Cref{thm:formalCDGA} 
is that $K \to G \to G/K \to BK \to BG$ is a fiber sequence and 
$G/K$ is the pullback of $BK \to BG \from EG$ up to homotopy,
allowing one to identify $H^*(G/K)$ with 
the cohomology of a \CDGA $H_K^* \ox H^*G = H_K^* \ox_{H_G^*} (H_G^* \ox H^*G)$, where $H_K^*$ has zero differential and 
$H_G^* \ox H^*G$ is a Koszul algebra modeling $EG$.\footnote{\ 
	This pullback argument is due to Borel's thesis~\cite[Thm.~25.2]{borelthesis},
	the model to Cartan~\cite[Thm.~5, p.~216]{cartan1950transgression}.
}

\bdefn
Let $\F$ be a field and $A$ a commutative unital $\F$-algebra.
A \defd{Noether normalization} of $A$ is a polynomial subring 
$B = \F[f_1,\ldots,f_n]$ of $A$ such that $B$ is a finite $A$-module.
Such a sequence $\vec f = (f_1,\ldots,f_n)$ is called
a \defd{system of parameters} for $A$.
Given a finite-dimensional representation $V$ of a finite group $\G$
over $\F$, the ring $\defm{\F[V]}$ of polynomial functions $V \lt \F$
naturally inherits a $\G$-action by ring automorphisms.
We grade $\F[V]$ by setting $\defm{|v^*|} = 1$ for $v^* \in \Hom_\F(V,\F)$
and write $\defm{|\G|}$ as well for the order of $\G$.
A \defd{pseudoreflection group} $\G \leq \Aut_\F V$ is a subgroup 
generated by elements fixing codimension-one subspaces.
\edefn


\bprop[{\cite[Prop.~5.5.5]{smithinvariantbook}}]\label{thm:invariantdegree}
Let $V$ be a finite-dimensional representation over a field $\F$
of a finite group $\G$.
If $\F[V]^\G$ contains a system of parameters $\vec f$
such that $\prod |f_j| = |\G|$, then $\F[V]^\G = \F[\vec f]$.
\eprop

\bprop[{\cite[p.~87 after Def.]{kane1994poincare}}]%
\label{thm:bound}
Let $V$ be a faithful finite-dimensional representation over a field $\F$
of a finite group $\G$.
Then $\dim_\F \big(\F[V]\quot \F[V]^\G\big) \geq |\G|$,
with equality if and only if the image of $\G$ is a pseudoreflection group.
\eprop

Recall by contrast that if $\vec f$ is a regular sequence of length $\dim V$ in $\F[V]$,
then 
\quation{\label{eq:CIRdim}
\dim_\F \F[V] \quot \F[\vec f] = \prod |f_j|.
}
This is because if $d$ is the left-hand side,
then $\F[V]$ is free of rank $d$ over $\F[\vec f]$.
The Poincar\'e series of $\F[V]$ and $\F[\vec f]$
are respectively $(1-t)^{-n} = (1+t+t^2+\cdots)^n$ and $\prod (1-t^{|f_j|})^{-1}$,
so $d$ is the 
value of the polynomial $\prod (1-t^{|f_j|}) / (1-t)^{-n}$ 
at $t = 1$.

\begin{lemma}\label{thm:shiga}
	Let $B$ be a connected \CGA over a field $\F$,
	admitting an action by algebra automorphisms by a finite group $\G$ of order
	relatively prime to the characteristic of $\F$,
	and let $A$ be a graded $\F$-subalgebra of the invariant ring $B^\G$.
	If $B\quot A =  B \quot B^\G$,
	then $A = B^\G$.
\end{lemma}
\bpf
Write $\defm{\mu}\: B \longepi B^\G$ 
for the $\G$-average $b \lmt \smash{\frac 1{|\G|} \sum_{\g \in \G} \g b}$.
%
As we assume the ideals of $B$ 
generated by $A^{\geq 1}$ and $\smash{\big(B^{\geq 1}\big)\mnn{}^\G}$
are equal,
any homogeneous element $b$ of positive degree in $\smash{B^\G}$ 
may be written as
$\sum a_j b_j $ for some $\smash{a_j \in A^{\geq 1}}$ and $b_j \in B$.
Since $f$ and the $a_j$ are invariant, 
\[
b = 
\mu(f) = 
\sum\, a_j \.\mu (b_j) .
\]
Thus $(B^\G)^{\geq 1} = A^{\geq 1} \. B^\G$.
By the graded Nakayama lemma~\cite[Prop.~5.2.3]{smithinvariantbook}, then, 
$A = B^\G$. 
%
\epf

\bthm[Chevalley--Shepherd--Todd \textup{\cite[Thm.~1.5 \& p. 82]{kane1994poincare}}]%
\label{thm:chevalleyshepherdtodd}
Let $V$ be a faithful representation of a finite group $\G$
over a field $\F$ of characteristic
relatively prime to the order of $\G$.
\TFAE:
\benum
\item The image of $\G$ is a pseudoreflection group.
\item The invariant ring $\F[V]^\G$ is polynomial.
\item The quotient $\F[V] \quot \F[V]^\G$ is a \CIR.
\item The quotient $\F[V] \quot \F[V]^\G$ is a Poincar\'e duality
algebra.
\eenum
\ethm

Now we have enough to prove the our first strengthening 
of the Shiga--Takahashi theorem.


\bpf[Proof of \Cref{thm:ST+}]
%
Since $G/K$ is formal
there is by \Cref{thm:formalCDGA} a regular sequence 
$\vec f$ in $H_K^*$
which lies in the image of $H_G^* \lt H_K^*$
and is such that $H_K^*\quot H_G^* = H_K^* \quot \Q[\vec f]$.
We have a sandwich
\[
\Q[\vec f] \leq \im(\mn H_G^* \to H_K^*) \leq \HKN.
\]
If $\HKN = \Q[\vec f]$, 
the sandwich collapses, 
showing $H_G^* \lt \HKN$ is surjective.
On the other hand, if $H_G^* \lt \HKN$ is surjective,
then $H_K^*\quot \HKN = H_K^*\quot H_G^* = H_K^* \quot \Q[\vec f]$,
so the sandwich collapses by \Cref{thm:shiga}.
In short $H_G^* \lt \HKN$ is surjective 
if and only if 
$\HKN = \Q[\vec f]$.

By \Cref{thm:invariantdegree}, if $|N| = \prod|f_j|$, then $\HKN = \Q[\vec f]$.
On the other hand, if $\HKN = \Q[\vec f]$,
then as $H_K^*\quot\HKN = H_K^* \quot \Q[\vec f]$
is a \CIR,
it follows by \Cref{thm:chevalleyshepherdtodd} that $N$ is a reflection group.
Then by \eqref{eq:CIRdim} and \Cref{thm:bound},
\[\prod |f_j| = \dim_\F H_K^* \quot \Q[\vec f] = \dim_\F H_K^* \quot \HKN = |N|.\]

Now, by \eqref{eq:CIRdim} and \Cref{thm:formalCDGA}, we know
\[
\prod |f_j| = 
\dim_\Q H_K^* \quot \Q[\vec f] = 
\dim_\Q  H_K^* \quot H_G^* = 
\frac{\dim_\Q H^*(G/K)}{2^{\rk G - \rk K}}
\]
since $G/K$ is formal.
Finally $H_G^* \lt \HKN$ is surjective 
if and only if $\HKN = \Q[\vec f]$,
if and only if $|N| = \prod |f_j| = \left.\dim_\Q H^*(G/K) \,\right/ {2^{\rk G - \rk K}}$;
but this last equation is equivalent to \isotfity by \Cref{thm:fpdim}.
\epf

With \Cref{thm:isotfformal} this immediately can be strengthened again:

\begin{restatable}{theorem}{thm:ST++}\label{thm:ST++}
	A \ccpair \GK is \isotf 
	if and only if 
	it is formal and $\rho^*\: H_G^* \lt \HKN$ is surjective.
\end{restatable}

In turn, using \ref{thm:chevalleyshepherdtodd} again,
one has a further enhancement:

\STPPP*
\begin{proof}
	The equivalence (1) $\iff$ (2) is \Cref{thm:ST++}.
	For (1) $\implies$ (3),
	apply \Cref{thm:isotfformal}
	and 
	the proof of \Cref{thm:ST+}.
	For (3) $\implies$ (1), 
	if $H_G^* \longepi \HKN$ and $N$ acts on $QH_K^*$ as a reflection group,
	then $H_K^* \quot H_G^* = H_K^* \quot \HKN$ is a \CIR, so by \Cref{thm:formalCDGA},
	the pair \GK is formal.
	To see (4) $\iff$ (1), recall from \Cref{thm:SGS}
	that $(G,K)$ is \isotf just if \GS is, and that there is a natural isomorphism
	from $H^2(BS;\R)$ to the dual space of $\fs$.
\end{proof}

Since \Cref{thm:SGS} and \Cref{thm:formaltorusreplacement}
let one exchange $S$ and $K$, 
one also finds the following invariant-theoretic results,
which do not seem to be otherwise obvious.

\bcor
Let $(G,K)$ be a formal pair and $S$ a maximal torus of $K$.
Then $H_G^* \lt H_K^{N_G(K)}$ is surjective if and only if
$H_G^* \lt H_S^{N_G(S)}$ is.
\ecor

\bcor
Let $(G,K)$ be a formal pair and $S$ a maximal torus of $K$.
Suppose $H_G^* \lt H_K^{N_G(K)}$ or $H_G^* \lt H_S^{N_G(S)}$ is surjective. 
Then $\pi_0 N_G(S)$ acts on the Lie algebra $\fs$ of $S$ as a reflection group
if and only if  
$\pi_0 N_G(K)$ acts as a reflection group on a (hence, any) 
minimal homogeneous vector space of algebra generators
of $H_K^*$.
\ecor

\bs

\nd\footnotesize{%
	\url{jeffrey.carlson@tufts.edu
	}
}
\end{document}